\documentclass[10pt]{article}

\usepackage{listings}
\usepackage[colorlinks=false, backref]{hyperref}

\usepackage{amsmath,amsthm,amsfonts,amssymb,amscd}
\usepackage{mathabx}
\usepackage{float}
\usepackage{graphicx}
\usepackage[usenames]{color}
\usepackage{makeidx}

\usepackage[mathscr]{euscript}
\usepackage{stmaryrd}
\usepackage{microtype}
\usepackage{booktabs}
\usepackage{cleveref}
\usepackage{bookmark}
\usepackage{mathrsfs}

\usepackage{emptypage}
\usepackage{tikz}

\DeclareMathAlphabet{\mathpzc}{OT1}{pzc}{m}{it}

\theoremstyle{plain}
\newtheorem{theorem}{\scshape Theorem}

\theoremstyle{definition}

\newtheorem{definition}{\scshape Definition}
\newtheorem{remark}{\scshape Remark}

\theoremstyle{definition}

\newcommand{\pset}[1]{ (#1) }

\newcommand{\norm}[1]{ \|#1\| }

\newcommand{\abs}[1]{ |#1| } 
\newcommand{\Abs}[1]{\left |#1\right| }

\newcommand{\cp}[1]{,_{#1}}

\def\hd{\bar{\partial}}

\def\d{\displaystyle}
\def\div{\operatorname{div}}

\def\n{\nonumber}

\def\p{\partial}

\def\e{\tilde\eta}
\def\v{\tilde v}
\def\a{\tilde a}
\def\A{\tilde A}
\def\J{\tilde J}
\def\os{{\Omega_0^s}}
\def\of{{\Omega_0^f}}
\def\g{\Gamma}
\def\u{\underbrace}

\def\ud{\frac{1}{2}}
\def\ddt{\frac{d}{dt}}

\usepackage{fancyhdr}
\title{Global existence and convergence near equilibrium for the moving interface problem between Navier-Stokes and the linear wave equation}

\author{
Daniel Coutand
\\Maxwell Institute for Mathematical Sciences and
\\Department of Mathematics, Heriot-Watt University,
\\Edinburgh, EH14 4AS  UK
\\{\footnotesize email: d.coutand@hw.ac.uk}
}

\begin{document}

\maketitle

\noindent
{\bf Abstract.} {\small We first establish existence for all positive time near equilibrium for the moving interface problem between the Navier-Stokes equations for the evolving fluid phase (moved by the fluid velocity) and an elastic body modelled by the linear wave equation. This problem has an infinite number of simple solutions with a flat interface (with zero velocity in the fluid, and zero horizontal velocity in the solid), that we call flat interface solutions. We then show that if the initial data is close enough to the canonical equilibrium, the solution converges towards a flat interface solution in large time, showing that these flat interface solutions capture the long time behaviour of this fluid-structure problem near the canonical equilibrium. This result is established with gravity (which can be set to zero or not).}


\section{Introduction}
 Fluid-structure interaction problems are omnipresent in nature, and involve at the mathematical level a coupling between a fluid phase (most often modelled by the Navier-Stokes equations) and a solid phase (most often modelled by an hyperbolic PDE), with a time dependent interface moving with the velocity of the fluid. 
Following the early works of \cite{Wein}, \cite{Serre}, \cite{EEM}, the analysis of the moving interface problem between a viscous fluid and a solid structure has become a quite active field of research.

The first type of problems studied is when the solid is rigid (and so does not change shape, but moves in and interacts with the fluid phase). For the motion of a rigid body inside a viscous fluid, existence of a weak solution (until eventual collision with the boundary) was established in \cite{DeEs1999}, and global in time existence of weak solutions (no longer unique if a collision occurs) was established in \cite{SMST} in $2$D and \cite{Fe2003} in $3$-D. Similar or related results were obtained independently by \cite{Conca}, \cite{FlOr}, 
\cite{Gunz}. Higher regularity properties for such models are obtained in \cite{Gr}, \cite{GrMa}.

  The impossibility of finite time collision with a fixed smooth boundary (under the Dirichlet boundary condition $u=0$ on the fixed boundary of the fluid domain), was established in \cite{Hi2007} and \cite{HiTa2009}. Surprisingly (in light of the previous result), if the rigid body has a boundary with limited regularity, \cite {GVHi2010} established that finite time contact is possible. Different boundary conditions are considered in \cite{GVHiWa2015} and \cite{GVHi2014}, and can also lead to different conclusions (\cite{GVHiWa2015}) than with the standard Dirichlet condition. The existence of zero acceleration configurations for the rigid body moving in the viscous fluid was established by \cite{HiSe}.

Global in (positive) time existence of weak solutions for fluid-elastic interaction has been obtained if the elastic model has a high enough degree of space regularity at the level of the basic energy law. Various models have been treated, with either finite dimensional restriction, or operators of degree at least four (with for instance Koiter type plates and shells). The models of elasticity can either be three dimensional, or plate/shell models: See \cite{DeEsGrLe}, \cite{ChDeEsGr}, \cite{LW}, \cite{Gunz}, \cite{FlOr}, \cite{BGLT1}, \cite{BGLT2}, \cite{B}, \cite{G}, \cite{LR}, \cite{BrSc}, \cite{MC1}, \cite{MC2}, \cite{BGRW}. These global results hold so long as there is no contact issue, such as self-intersection of a moving boundary, or intersection between a moving boundary and a fixed boundary. In \cite{GH}, it is shown further that a damped version of the beam model considered in \cite{MC1} (for which global in time existence of weak solutions, assuming no contact issues,  follows similarly as in \cite{MC1}) actually does not have any contact issue developing in time. The method needs a strong solution approach, with more regularity on the initial data, and a two-dimensional setting for Navier-Stokes.

The most canonical model for a three dimensional deformable elastic phase is either the linear wave equation, or the linear system of elasticity with Lam\'e coefficients, which are second order hyperbolic PDE and do not feature the strong a priori control offered by fourth order operators described in the previous paragraph. For this reason, weak solutions for the interaction between this  classical model for the solid phase coupled with Navier-Stokes for the moving fluid phase are not known to exist.
 
The existence of local in time solutions to the interaction between the incompressible Navier-Stokes equations and linear or nonlinear elasticity was first established in [\cite{CS1}, \cite{CS2}]. Further results, with more general initial data, were established in \cite{KT1}, \cite{KT2}, \cite{KTZ}, \cite{RV}, \cite{BGT}. Local in time existence for the interaction between the
compressible Navier-Stokes equations and elasticity was first established in \cite{BG1}, \cite{BG2}, \cite{KT3}. See also \cite{BDV}, \cite{ALT}, \cite{AT1} and \cite{AL} for other existence results of strong solutions in fluid-elastic interaction problems.

Global existence (near equilibrium) for this type of moving interface problems coupling a Navier-Stokes fluid phase with a second order hyperbolic equation in the solid phase, was still completely open until the present paper. The results previously obtained so far in this direction have considered models with damping being added either in the solid phase or on the interface, see for instance \cite{IKLT2}, \cite{IKLT3}. The recent results \cite{KO1}, \cite{KO2} treat the case of the lightest and most natural damping being added to the linear wave equation, by adding a friction type term $-\alpha v$ as forcing (this can be thought of as an air friction term being applied in the elastic phase to add dissipation). The addition of damping to the linear wave equation in the solid phase provides directly some dissipative terms in this phase, although some delicate analysis is needed to obtain the remaining relevant quantities. Without this addition, there is no obvious way to identify dissipative terms in the solid phase, which renders the issue of global existence substantially more challenging than the already complex cases previously treated in the literature.

The aim of the present paper is first to establish global in (positive) time existence for the moving interface problem between the incompressible Navier-Stokes equations and the linear wave equation modelling the elastic phase (without any damping being added), when the initial data is close to an equilibrium.

 We then show convergence in large time toward a particular solution of the problem with a flat interface (there is an infinite number of such particular solutions). This is done with the gravity constant $g$  which is either strictly positive or null. The paper establishes these results when the volume of the elastic phase is close enough to the volume of its reference configuration (where the linear wave equation is set).

We now briefly outline the paper.

In Section \ref{sec:notation}, we define our notations. In particular, we take the viewpoint of defining fractional Sobolev norms on $\Gamma$ via Fourier series, which is very useful later on, in Section \ref{crucial}.

In Section \ref{description}, we introduce the problem, which is set with {\bf periodic boundary conditions in the canonical horizontal directions}. We then identify the canonical equilibrium. This problem also has an infinite number of special solutions, with zero velocity in the fluid phase, a flat interface, and where the dynamics in the solid is dictated by a one dimensional linear wave equation with homogenous Dirichlet conditions for the vertical component of the displacement. We call these special solutions {\it flat interface solutions} of the problem, and we establish later in this paper that they describe the behaviour of the system in large time, when the initial data is close enough to the canonical equilibrium.

In Section \ref{ALE}, we introduce our Arbitrary Lagrangian representation of the problem via a Stokes extension problem (and not the standard Lagrangian coordinates in the fluid), in order to obtain  global in time existence. The fundamental interest of this choice of representation is explained in Section \ref{crucial}. We then introduce our functional framework, with in particular a dissipative norm (associated with the viscosity in the fluid). We finally state our main Theorems and  write down the variational problems respectively associated with the system (\ref{ale}) and its first and second time differentiated versions. 

In Section \ref{sectionpressure}, we gather needed estimates for the pressure that will be used later on.

In the crucial Section \ref{crucial}, we show that first order horizontal derivatives of the displacement $\eta$ (evaluated from the solid phase $\os$) are controlled on the boundary in $L^2(0,T;H^{\frac 3 2}(\g))$ in terms of the dissipative norm in the fluid. This is a surprising result, as $\eta$ is not part of the dissipative norm a priori. This is why we choose to work with an extension operator of $\eta$ from $\g$ to $\of$ with an Arbitrary Lagrangian representation, as a priori the standard Lagrangian map in the fluid may not have this dissipative behaviour.

In Section \ref{highspace}, we perform the highest in space order estimates on the zero time differentiated problem. 

In Section \ref{middle}, we explain why the first time differentiated problem is the one where no issue arises, and state the estimate associated.

In Section \ref{end}, we then perform the estimates associated with the highest in time problem. We also establish and state higher regularity properties being satisfied by the solution of the problem, which although not necessary, make estimates lighter.

We conclude the proof of Theorem \ref{theorem_main} (namely the existence of the solution for all positive time if the initial data is close enough to the canonical equilibrium) in Section \ref{conclusion1}.

For a solution satisfying the assumptions of Theorem \ref{theorem_main}, we then establish Theorem \ref{theorem_asymptotic}, namely the convergence of the solution in large time towards a {\it flat interface solution}, in Section \ref{conclusion_asymptotic}. We first show that the solution in the fluid phase converges towards zero, and that the interface converges towards the flat interface. We then study how these properties propagate in the solid phase to some extent, by continuity of velocity and normal stress at the interface. The limit problem for the vertical displacement in the solid phase is identified as a weak limit of a sequence of one dimensional wave equations with homogenous Dirichlet boundary conditions and defined at time $n$ with the horizontal average of the solution to Theorem \ref{theorem_main}.

\section{Notations} \label{sec:notation} 

\subsection{Notation for the gradient vector} \label{sec:grad-horiz-deriv}

Throughout the paper the symbol $\nabla $ will be used to denote the three dimensional gradient vector 
$
\nabla =\left( \frac{\p}{\p x_1}\,,  \frac{\p}{\p x_2}\,,  \frac{\p}{\p x_3}   \right)
$.
\subsection{Notation for partial differentiation and  the Einstein summation convention} \label{sec:notat-part-diff}

The first time derivative of $F$ will be denote by $F_t=\frac{\partial  F}{\partial t}$, while the second and third time derivatives will be denoted $F_{tt}$ and $F_{ttt}$.

 The $k$th spatial partial derivative of $F$ will be denoted by $F\cp{k}=\frac{
	\partial F}{
	\partial x_k}$. Repeated Latin indices $i,j,k$, etc., are summed from $1$ to $3$.

For example, $\d F\cp{ii}=\sum_{i=1}^3\frac{\p^2F}{\p x_i\p x_i}=\Delta F$, and $\d F^i\cp{j} I^{jk} G^i\cp{k}=\sum_{i=1}^3\sum_{j=1}^{3}\sum_{k=1}^3\frac{\p F^i}{\p x_j} I^{jk} \frac{\p G^i}{\p x_k}$.

\def\R{ \mathbb{R}  }

\subsection{Tangential (or horizontal) derivatives}\label{sec: tangential derivative} 

Depending on context, we define $\hd$ as either being simply an horizontal derivative $\hd = \frac{\partial}{\p x_\alpha}$, $\alpha=1,2$, or as being the horizontal gradient $\hd=\left( \frac{\p}{\p x_1}\,,  \frac{\p}{\p x_2} \right)$.

\subsection{Sobolev spaces} \label{sec:diff-norms-open}

For integers $k\ge0$ and a bounded set $U$ of $\R^3$ under the form $U=(0,L)^2\times [h_1,h_2]$, we define the Sobolev space $H^k(U)$ $\pset{H^k(U;\R^3)}$ to be the completion of the set of functions {\bf periodic with period $L>0$ in the canonical horizontal directions} $d_1=(1,0,0)$ and $d_2=(0,1,0)$, and of regularity  $C^\infty(\bar{U})$ $\pset{C^\infty(\bar{U}; \mathbb{R}  ^3)}$ in the norm 
\begin{align*}
	\norm{u}_{H^k(U)}^2=\sum_{\abs{a}\le k}\int_U \Abs{ \nabla ^a u(x) }^2 dx=\sum_{\abs{a}\le k}\int_U \Abs{ \left(\frac{\p}{\p x_1}\right)^{a_1}\left(\frac{\p}{\p x_2}\right)^{a_2}\left(\frac{\p}{\p x_3}\right)^{a_3} u (x) }^2 dx  , 
\end{align*}
for a multi-index $a\in \mathbb{Z}  ^3_+$, with the convention that $\abs{a}=a_1+a_2+a_3$. 
For conciseness, we will write $H^s(U)$ instead of $H^s(U;\R^3)$ for vector-valued functions. 

\begin{remark}
	All our functions will be periodic of period $L$ in the canonical horizontal directions $d_1$ and $d_2$, so this is why for conciseness of notations we choose not to add any additional index symbol for periodicity (such as $H^k_\text{per}(U)$ for instance) for our definition of $H^k(U)$ in this paper.
\end{remark}

\subsection{Sobolev spaces on $\Gamma$, or on any horizontal plane} \label{sec:sobolev-spaces-gamma}

We denote by $e_n$ the Hilbert basis of orthonormal functions of the variable $x_h=(x_1,x_2)$, $L-$periodic in $x_1$ and $x_2$, and with zero average, made of eigenfunctions of the two-dimensional Laplacian with periodic boundary conditions. Namely, each $e_n$ is a product of sine or cosine functions periodic with period $L$ in $x_1$, $x_2$ and satisfies
$$ \Delta_0{e_n}=e_n,_{11}+e_n,_{22}=-\lambda_n e_n\,.$$

Any smooth $L-$periodic in $x_1$ and $x_2$ function $f$ with zero average on $\g$ (or on any horizontal slice $[0,L]^2\times \{x_3\}$) can be expanded in this basis:
\begin{equation*}
	f(x_h,h_s)=\sum_{n=1}^\infty f^n e_n(x_h)\,,
\end{equation*}
and the $H^s$ norm of $f$ on $\g$ ($s\ge 0$) can be classically defined as:
$$\|f\|_{H^s(\Gamma)}^2=\sum_{n=1}^\infty \lambda_n^{s} (f^n)^2\,.$$
The typical example of function with zero average we will use are horizontal derivatives (of periodic functions).

If $f$ is not of zero average on $\g$, we define for $s\ge 0$
$$\|f\|_{H^s(\Gamma)}^2=\left|\frac{1}{|\Gamma|}\int_\g f dx_h\right|^2+\sum_{n=1}^\infty \lambda_n^{s} (f^n)^2\,.$$

\subsection{The unit normal vectors}

If we work in the solid reference domain $\os$ described in the next Section, we will denote by $N=N^s=(0,0,1)$ the outward unit normal vector.

If we work in the fluid reference domain $\of$ described in the next Section, we will denote by $N=N^f=(0,0,-1)$ the outward unit normal vector.

We will only add a superscript $N^s$ or $N^f$ when needed to avoid confusion.

If we work in the moving domain $\Omega^s(t)$ or $\Omega^f(t)$, the same type of convention is adopted for the normal outward unit vector $n$, which is either $$n^s=\frac{(a_i^3)_{i=1}^3}{|(a_i^3)_{i=1}^3|}=\frac{X,_1\times X,_2}{|X,_1\times X,_2|}$$ or $n^f=-n^s$, where $X$ and $a$ are defined in the next Section.

\subsection{Conventions on constants and polynomials}\label{polynomials}

We take the convention that $P$ denotes a generic polynomial of degree greater or equal than $1$, with positive coefficients which do not depend on time or on how large $\lambda>0$ (elasticity coefficient in (\ref{307.1})) is. Moreover, the coefficient of the power zero for our generic $P$ is always equal to zero.

In the same way, $P_{\lambda}$ denote a polynomial similar as $P$, but with coefficients dependent on $\lambda$ (and getting possibly large for $\lambda$ large).

The same convention is adopted for positive constants $C$ and $C_\lambda$.

The necessity to keep track of the dependence in $\lambda$ is for the case of non zero gravity, as we will see $\lambda$ needs to be large enough relative to $g$.

For estimates where this tracking is not needed, we will just use the standard $A\lesssim B$ to indicate there exists $c>0$ finite independent of time such that $A\le c B$.

\subsection{Cofactor matrix and Piola's identity} \label{cofactor}

For a bijective mapping  $\tilde X =\text{Id}+\e$ (close to identity in this paper), we define $$\A=(\nabla \tilde X)^{-1}=\nabla(\e+\text{Id})^{-1}(\tilde\eta+\text{Id})\,,\hskip 0.3 cm\a=\text{Cof}\nabla(\e+\text{Id})\,,\hskip 0.3cm \J=\text{det}\nabla(\e+\text{Id})\,,$$ so that
$\d\A=\frac{\a}{\J}\,.$

	One crucial relation we will use later on will be { Piola's identity for the Cofactor matrix} (valid for each fixed $i$):
	\begin{equation}
		\a_i^j,_j=0\,.\label{piola}
	\end{equation}

\section{Description of the Navier-Stokes-linear wave interaction problem, special solutions with a flat interface}
\label{description}
\subsection{description of the problem}

\begin{figure}[h]
	\begin{tikzpicture}[scale=0.19]
		
		\draw[color=green,ultra thick] plot[smooth,tension=0.5] coordinates{( -1,-1)(2,-1.3)  (7,-0.7)(10,-0.7) (13, -1.3) (16, -0.7) (19,-1.3)(21,-0.7) (23,-1.3) };
		\draw[black,dashed] plot[smooth,tension=0.5] coordinates{(0,-10)(-0.2,-5)( -1,0) (-0.2,5) (0,10) };
				\draw[black,dashed] plot[smooth,tension=0.5] coordinates{(24,-10)(23.8,-5)( 23,0) (23.8,5) (24,10) };
		\draw[color=red,ultra thick] (0,10)-- (24,10);
		\draw[color=red,ultra thick] (0,-10)-- (24,-10);
		\draw (14,-6) node { { \small $\Omega^s(t)=\{X(x,t); x\in\Omega_0^s\}$}}; 
		\draw (14,6) node { { \small $\Omega^f(t)=\{X(x,t); x\in\Omega_0^f\}$}}; 
		\draw (13,-2.5) node { {\small $\Gamma(t)=\partial\Omega^s(t)\cap\partial\Omega^f(t)$}};
		\draw (14,11) node { {\small $\Gamma_{top}=\{x_3=h\}$}};  
		\draw (10,-11) node { {\small $\Gamma_B=\{x_3=0\}$}};
		\draw (0,-11) node { {\small $x_1=0$}}; 
		\draw (22,-11) node { {\small $x_1=L$}};
		\draw[thick,->] (2,-1.3) -- (1.8,-4);
		\draw (2.5,-2) node { { $n$}}; 
		\draw (12,-13) node { {\small Physical domain at time $t\ge 0$}};
		\draw (12,-15) node { {\small Elastic phase of average height $h_e$}};
		\draw[color=blue,ultra thick] (26,0)-- (50,0) ;
		\draw[black,dashed] (26,-10)-- (26,10);
		\draw[black,dashed] (50,-10)-- (50,10);
		\draw[color=red,ultra thick] (26,10)-- (50,10);
		\draw[color=red,ultra thick] (26,-10)-- (50,-10);
		\draw (40,-6) node { { \small $\Omega_0^s$}}; 
		\draw (40,6) node { {\small $\Omega_0^f$}}; 
		\draw (40,-1) node {{\small $\Gamma=\{x_3=h_s\}$}};
		\draw (40,11) node { { \small $\Gamma_{top}=\{x_3=h\}$}};  
		\draw (38,-11) node { {\small $\Gamma_B=\{x_3=0\}$}};
		\draw (28,-11) node { {\small $x_1=0$}}; 
		\draw (48,-11) node { {\small $x_1=L$}};
		\draw[thick,->] (32,0) -- (32,-3);
		\draw (40,-3) node { {\small $N=(0,0,-1)=-d_3$}}; 
		\draw (39,-13) node { {\small Reference domain}}; 
		\draw (39,-15) node { {\small Reference solid domain of height $h_s$}};
		\draw[color=green,ultra thick] (52,-1)-- (76,-1) ;
		\draw[black,dashed] (52,-10)-- (52,10);
		\draw[black,dashed] (76,-10)-- (76,10);
		\draw[color=red,ultra thick] (52,10)-- (76,10);
		\draw[color=red,ultra thick] (52,-10)-- (76,-10);
		\draw (65,-6) node { {\small $\Omega^s(\infty)=(0,L)^2\times (0,h_e)$}}; 
		\draw (65,6) node { {\small $\Omega^f(\infty)=(0,L)^2\times (h_e,h)$ }}; 
		\draw (65,-2) node {{\small $\Gamma_e=\{x_3=h_e\}$}};
		\draw (65,11) node { {\small $\Gamma_{top}=\{x_3=h\}$}};  
		\draw (63,-11) node { {\small $\Gamma_B=\{x_3=0\}$}};
		\draw (54,-11) node { {\small $x_1=0$}}; 
		\draw (74,-11) node { {\small $x_1=L$}};
		\draw[thick,->] (32,0) -- (32,-3);
		\draw (65,-13) node { {\small Asymptotic physical domain}}; 
		\draw (65,-15) node { {\small Asymptotic solid of height $h_e$}};
	\end{tikzpicture} 
\end{figure}

We consider a fixed domain $\Omega=(0,L)^2\times (0,h)$ ($L>0, h>0$), with {\bf periodic boundary conditions of period $L$ in the horizontal directions} $d_1=(1,0,0)$ and $d_2=(0,1,0)$.

 This domain is partitioned in two reference subdomains $\os=(0,L)^2\times (0,h_s)$, with $h_s\in (0,h)$ (for the solid phase) and $\of=(0,L)^2\times (h_s,h)$ (for the fluid phase), with a flat reference interface $\Gamma$ at $x_3=h_s$.
 \begin{remark}
  In the context of the reference domains, the subscripts $0$ do not mean the initial configuration of the solid and fluid phases have to be exactly at $\os$ and $\of$, although in this paper we will assume they start close to this configuration. We will then show that for initial data close enough to the canonical equilibrium, the solution of the problem is defined for all positive time, and that the geometry of the fluid and solid phases converge towards the right picture above.
 \end{remark}
 
 \begin{remark}
 		The volume of the solid phase remains constant for all time due to incompressibility in the fluid.
 	
 \end{remark}

  The elastic phase $\Omega^s(0)$ at time $t=0$ is described by the {\bf position} map $X(x,0)=x+\eta(x,0)$  ($x\in\os$). In this paper, the {\bf displacement} $\eta$ will be assumed small initially. At time $t\ge0$, the solid domain (from the basic periodic cell) is $\Omega^s(t)$. A particle initially at $x+\eta(x,0)$ is located at time $t$ at the position $X(x,t)=x+\eta(x,t)$  ($x\in\os$). Our linear elasticity model is assumed to satisfy the classical linear wave equation:
  \begin{equation}
  	v_t-\lambda \Delta \eta=-g d_3\,\ \ \text{in}\ \os\times [0,T]\,,\label{307.1}
  \end{equation}
  where $g\ge 0$ is the gravity constant (null, if no gravity is considered), $d_3=(0,0,1)$, $\lambda>0$ and the velocity $v=\eta_t$. This PDE is complemented with the Dirichlet boundary condition on $\Gamma_B=[0,L]^2\times \{0\}$:
  \begin{equation}
  	\eta=0\,\ \ \text{on}\ \Gamma_B\times [0,T]\,,\label{307.2}
  \end{equation}
  and an initial condition
  \begin{align}
  \eta(\cdot,0)=&\eta_0(\cdot)\,,\\
  v(\cdot,0)=&v_0(\cdot)\,,
  \end{align}
  and further boundary conditions on the interface $\Gamma$ that we will describe later in this Section.
 
 The average height of the initial solid domain is $h_e\in (0,h)$ given by 
 \begin{equation}
 	\label{averageheight}
 	h_e L^2=|\Omega^s(0)|\,.
 \end{equation}
 
 We also denote 
 \begin{equation}
 	\label{equilibriumboundary}
 	\Gamma_e=[0,L]^2\times \{h_e\}. 
 \end{equation}
 
 \begin{remark}
 	The case when $h_e<h_s$ (respectively $h_e>h_s$) corresponds to the case when the initial solid domain is on average compressed (respectively extended) compared to its reference configuration $\os$, and the case $h_e=h_s$ corresponds to the case when the initial solid domain has the same volume as the reference configuration $\os$. Due to the conservation of volume in the fluid phase (incompressibility), the solid phase keeps the same volume for all time, and the same average height.
 \end{remark}

 \begin{remark}
 Note that the average height of the domain $\Omega^s(0)$ is $h_e$, and is assumed close to the height $h_s$ of the reference domain where the wave equation is set.
 \end{remark}

The fluid phase $\Omega^f(0)$ at time $t=0$ is described by the position map $x+\eta(x,0)$  ($x\in\of$). In this paper, $\eta$ will be assumed small initially. A particle initially at $X(x,0)=x+\eta(x,0)$ is located at time $t$ at the position $X(x,t)=x+\eta(x,t)$  ($x\in\of$).  

For $0 \le t \le T$,
the evolution of the incompressible viscous fluid in the moving phase $$\Omega^f(t)=X(\Omega^f_0,t)$$
is modeled by the in\-com\-pres\-sib\-le Navier-Stokes equations, together with the boundary and initial conditions:
\begin{subequations}
  \label{NSe}
\begin{alignat}{2}
u_t+ u\cdot \nabla u -\nu\Delta u +  \nabla  p&= -g d_3  \ \  \ &&\text{in} \ \ \Omega^f(t) \,,\label{NSe.a}\\
  {\operatorname{div}} u &=0 
&&\text{in} \ \ \Omega^f(t) \,, \label{NSe.b}\\
u   &= 0  \ \  &&\text{on} \ \ \Gamma_{top}=\{x_3=h\} \,,\label{NSe.c}\\
\nu \nabla_n u - p\, n &= -\lambda \eta,_3^s (X^{-1}) \ \ &&\text{on} \ \ \Gamma(t) \,, \label{NSe.d}\\
u   &= v^s(X^{-1})  \ \  &&\text{on} \ \ \Gamma(t) \,,\label{NSe.e}\\
u   &= u_0  \ \  &&\text{on} \ \ \Omega^f(0) \,,\label{NSe.f}
\end{alignat}
\end{subequations}
where $\nu>0$,  $v^s$ in (\ref{NSe.e}) denote the trace of the velocity $v$ defined in $\os$, and $\eta,_3^s$ in (\ref{NSe.d}) denote the trace of $\eta,_3$ defined from $\os$. Generally, in this paper we will denote velocity fields without superscripts relating to their phase, except when the context requires precision.
We precise that 
 $n(\cdot,t)$ denotes the exterior unit normal vector to $\Omega^f(t)$ at the moving interface  $\Gamma(t)=\p\Omega^f(t)\cap\p\Omega^s(t)$.
  The vector-field $u = (u_1, u_2, u_3)$ denotes the Eulerian velocity
field, and $p$ denotes the pressure function.

\begin{remark} 
	To make notations lighter, we use the standard Neumann boundary condition in (\ref{NSe.d}), instead of the symmetric gradient. This does not make the problem different, or easier mathematically.
	\end{remark}
	
\begin{remark}
Note that we take the viewpoint of setting identical density (to $1$) in the fluid and solid phases, mainly to have  more condensed notations later on. This restriction does not matter, and the proof presented here carries for the case of different phase densities.
\end{remark}
Here, (\ref{NSe.a}) and (\ref{NSe.b}) are the incompressible Navier-Stokes equations, written in Eulerian variables in the moving domain $\Omega^f(t)$. We assume a Dirichlet boundary condition (\ref{NSe.c}) on the flat top of the domain. The condition (\ref{NSe.e}) is stating the continuity of velocity across the interface, whereas (\ref{NSe.d}) states the continuity of normal stress across the moving interface. The final condition of (\ref{NSe}) are about the initial Eulerian velocity in the initial fluid phase.

The pressure $p$  is a solution to the following Dirichlet and Neumann problem:
\begin{subequations}
  \label{p}
\begin{alignat}{2}
- \Delta p  &=  u^i,_j u^j,_i   \ \  \ &&\text{in} \ \ \Omega^f(t) \,,\\
 p &=  n \cdot [\nu\nabla_n u+{\lambda} \eta,_3^s(X^{-1})] \ \ &&\text{on} \ \ \Gamma(t)  \,,\\
 \frac{\partial p}{\partial x_3}&=\nu \Delta u^3-g\ \ &&\text{on} \ \ \Gamma_{top}  \,,
\end{alignat}
\end{subequations}
so that given an initial domain $\Omega^f(0)$ and an initial velocity field $u_0$, the initial pressure is obtained as the solution of (\ref{p}) at $t=0$.

As we can see from our equations, the solid phase is naturally cast in the reference domain $\os$, whereas the fluid phase in its most usual form is cast in the moving domain $\Omega^f(t)$. It is more convenient to write all equations in the fixed reference domains. The most natural way to do so is usually to introduce the Lagrangian variables associated to the fluid phase (although we will not carry our analysis in this representation).

If we denote by $X=\text{Id}+\eta$ the flow map associated with $u$ defined in $\Omega^f(t)$, we have for all $x\in\of$, 
\begin{align*}
X_t=&u(X(x,t),t)\,,	\\
X(x,0)=&x+\eta(x,0)\,.
\end{align*}

Note that we can choose a different map $\tilde X$ representing the fluid domain $X(\R^2\times\{h_s,h\},t)$ from $\R^2\times\{h_s,h\}$, instead of the Lagrangian flow map $X$, so long as the map coincides with $X$ on the boundary of $\R^2\times\{h_s,h\}$. The main difference would be that it adds an extra advection term, which is a small perturbation in the small data regime. Another difference is that for the elementary cell, $X(\of,t)$ and $\tilde X(\of,t)$ would not necessarily be equal, but this does not matter by periodicity, and by the fact that $X(\g,t)=\tilde X(\g,t)$. This is the method of Arbitrary Lagrangian representation. This is indeed the viewpoint of this paper, as it is needed for some estimates to represent the moving fluid domain from a Stokes extension of the displacement from the solid phase.

For such a mapping $\tilde X =\text{Id}+\e$ satisfying $\e=\eta^s$ on $\Gamma$ and $\e=0$ on $\Gamma_{top}$, we define $$\A=(\nabla \tilde X)^{-1}=\nabla(\e+\text{Id})^{-1}(\tilde\eta+\text{Id})\,,\hskip 0.3 cm\a=\text{Cof}\nabla(\e+\text{Id})\,,\hskip 0.3cm \J=\text{det}\nabla(\e+\text{Id})\,,$$ so that
$\d\A=\frac{\a}{\J}\,.$ We next define the Arbitrary Lagrangian velocity field and pressure in the reference domain $\of$ by $$v=u\circ(\e+\text{Id})\,,\hskip 1cm q=p\circ(\e+\text{Id})\,.$$

Since $\eta^s=\tilde\eta$ on $\g$, we then have from (\ref{NSe.e}) continuity of velocity $$v=v^s$$ across $\g$.

Due to this continuity, we will generally denote the velocity simply $v\in H^1(\Omega)$ (and just precise $v^f$ or $v^s$ when needed, for instance when a normal derivative is involved on $\g$).

The divergence free condition then becomes in $\of\times [0,T]$:
\begin{equation}
	\label{divale}
	\A_i^j v,_j^i=0=\a_i^j v,_j^i\,.
\end{equation}

The problem for $v(\cdot,t)\in H^1(\Omega)$ then becomes in Arbitrary Lagrangian variables:

\begin{subequations}
	\label{ale}
	\begin{alignat}{2}
		\J v_t^i+(v-\tilde v)_l \a_l^j v,_j^i-\nu \a_j^l (\A_j^k v^i,_k),_l+\a_i^j q,_j&= -\J g \delta_3^i \ \  \ &&\text{in} \ \ \of\times [0,T] \,,\label{ale.a}\\
		\A_i^j v,_j^i &= 0  \ \  &&\text{in} \ \ \of\times [0,T] \,,\\
		v &= 0\ \  && \text{on} \ \ \Gamma_{top}\cup\Gamma_{B}\times [0,T] \,,\label{ale.c}\\
		-\nu \a_j^l \A_j^k {v^f},_k^i N_l+\a_i^j q N_j&= - \lambda N_l \eta,_l^s\ \  &&\text{on} \ \ \g\times [0,T] \,,\label{ale.d}\\
		v_t-\lambda\Delta\eta&=-g d_3  \ \  &&\text{in} \ \ \os\times [0,T] \,,\\
		(v(\cdot,0), \eta(\cdot,0))&=(v_0(\cdot), \eta_0(\cdot))\ \ &&\text{in} \ \ \Omega\times\os \,.
	\end{alignat}
\end{subequations}  

Note that our requirement $v(\cdot,t)\in H^1(\Omega)$ ensures that the continuity (\ref{NSe.e}) of velocity fields across the interface is satisfied, and for this reason, we do not need to report it in (\ref{ale}). So, this is a case when a boundary condition is encoded in the definition of the functional framework. We could also have removed (\ref{ale.c}), and precised that $v(\cdot,t)\in H^1_0(\Omega)$.

\begin{remark} The precise form of $\e$ will be given in the elliptic system (\ref{he}).\end{remark}

\subsection{The canonical equilibrium of height $h_s$}

It is easy to see that $v=0$, $\e=0$ in $\Omega\times [0,\infty)$, together with 
\begin{subequations}
		\label{0911}
	\begin{align}
\eta=&\eta_e=(0,0,\frac{g}{2\lambda} x_3(x_3-h_s))\,,\ &&\text{in}\ \os\,,
\label{0911.1}\\
\eta=&0\,,\ &&\text{in}\ \of\,,
\end{align}
\end{subequations}
and 
\begin{equation}
q=q_e=-g(x_3-h_s)-\lambda \frac{g}{2\lambda} h_s=-g x_3+\frac{g}{2} h_s\,,\ \text{in}\ \of\,,\label{0911.2}
\end{equation}
is an equilibrium satisfying the system (\ref{ale}) for the initial data $(v(\cdot,0),\eta(\cdot,0))=(0,\eta_e(\cdot))$. Also note that
\begin{equation}
	q_e=-\ud g h_s\ \text{on}\ \g\,.\label{1407.1}
\end{equation}

Given that the canonical equilibrium satisfies on $\g$
\begin{equation*}
	q_e \delta_i^3 =-\lambda \eta_e,_3^i\,,
\end{equation*}  
and $N=(0,0,-1)$, a useful way (later on) to rewrite the continuity of stress across the reference interface (\ref{ale.d}) is first
\begin{equation}
	\label{systemhorizontal}
	-\nu \a_j^3 \A_j^k {v^f},_k^i +\a_i^3 q = - \lambda \eta^s,_3^i\,,
\end{equation}
which then becomes
\begin{equation}
	-\nu \a_j^3 \A_j^k {v^f},_k^i +\a_i^3 (q-q_e) +\a_i^3 q_e- \delta_i^3 q_e = - \lambda (\eta^s-\eta_e),_3^i\,. \label{0211.1}
\end{equation} 
 \begin{remark}
 Note this equilibrium is not possible to achieve for the case when the average height $h_e$ of the elastic phase is distinct from $h_s$. It is nevertheless convenient to use later in the analysis (see for instance (\ref{0711.1}) and (\ref{0711.4})), as it satisfies $\lambda\Delta\eta_e^3=g$ and $\eta_e^3=0$ on $\Gamma$.
 \end{remark}

A general family of particular solutions of (\ref{ale}) with a flat interface can be described as follows.
\subsection{A general class of special solutions with flat interface of height $h_e$}\label{flat}
With our definition (\ref{averageheight}), we still have $v=0$, but this time $\d\e=(0,0,(h_e-h_s) \frac{x_3-h}{h_s-h})$ in the stationary fluid phase $\of$ for all time, so that

 $\nabla \tilde X=\begin{pmatrix}
	1 & 0 &0 \\
	0 & 1 & 0 \\
	0 & 0& \d\frac{h_e-h}{h_s-h}\\ 
\end{pmatrix}$, $\tilde J= \d\frac{h_e-h}{h_s-h}$, $\A=\begin{pmatrix}
1 & 0 &0 \\
0 & 1 & 0 \\
0 & 0& \d\frac{h_s-h}{h_e-h}\\ 
\end{pmatrix}\,, $  $\a=\begin{pmatrix}
\d\frac{h_e-h}{h_s-h} & 0 &0 \\
0 & \d\frac{h_e-h}{h_s-h} & 0 \\
0 & 0& 1\\ 
\end{pmatrix}\,. $

 This time, in the elastic phase $\os$, $\eta$ is under the form

$$\eta=(0,0,\alpha(x_3,t))$$ with $\alpha$ solution of a one dimensional wave equation (in the vertical direction):
\begin{subequations}
	\label{wave1d}
	\begin{alignat}{2}
		\alpha_{tt}-\lambda\alpha,_{33} &=-g \ \  \ &&\text{in} \ \ (0,h_s)\times [0,\infty) \,,\\
		\alpha &=h_e-h_s \ \  &&\text{on} \ \ \{h_s\}\times [0,\infty) \,,\\
		\alpha &=0 \ \  &&\text{on} \ \ \{0\}\times [0,\infty) \,,\\
		\alpha(0) &=\alpha_0 \ \  && \text{on} \ \ (0,h_s)  \,,\\
		\alpha_t(0) &=\alpha_1 \ \  && \text{on} \ \ (0,h_s)  \,,
	\end{alignat}
\end{subequations}
where $\alpha_0$ and $\alpha_1$ are given initial data satisfying compatibility conditions $\alpha_0(h_s)=h_e-h_s$, $\alpha_0(0)=0$ and $\alpha_1=0$ on $\{0\}\cup\{h_s\}$.

The pressure in the fluid domain is then given by
\begin{equation}
	q(x,t)=-g \d\frac{h_e-h}{h_s-h} (x_3-h_s) -\lambda \alpha,_3(h_s,t) \,,\ \text{in}\ \of\,.
\end{equation}

We call these special solutions {\bf flat interface solutions} of (\ref{ale}), for the initial data $v(x,0)=1_{\Omega_s}(x)\alpha_1(x_3)$ and  $\eta(x,0)=\alpha_0(x_3)$. 

\begin{remark}
For notational convenience, we will also define for $x=(x_1,x_2,x_3)$, $\alpha(x)=\alpha(x_3)$, $\alpha_i(x)=\alpha_i(x_3)$ ($i=0,1$) when $x_3\in [0,h_s]$
\end{remark}

\begin{remark}
The canonical equilibrium is a particular case with $\alpha_0=\eta_e^3$, $\alpha_1=0$ and $h_e=h_s$.
\end{remark}

We will first show that if a solution starts close enough to the equilibrium (\ref{0911}), (\ref{0911.2}), it is defined over $[0,\infty)$, and remains close to it. We will then show that as $t\rightarrow\infty$, the solution converges (in a weaker norm than the norm of existence) towards a flat interface solution of (\ref{ale}), as defined in Section \ref{flat}, showing that these flat interface solutions describe the long time behaviour of the small data regime.

\section{Arbitrary Lagrangian formulation of the problem and statement of our results}\label{ALE}

The main difficulty that the standard Lagrangian formulation of the problem would introduce is that the obvious identity
$$\eta(x,t)=\eta(x,0)+\int_0^t u(X(x,s),s) ds\,,$$ would yield by Cauchy-Schwarz a multiplicative constant of order $\sqrt{t}$ when estimating the cofactor matrix for the Lagrangian variables, which is unsuitable for large time.

This is why we use an Arbitrary Lagrangian formulation of the problem, where the fluid domain is described not by using the Lagrangian coordinates, but instead by using the Stokes  extension in $\of$ of the trace of $\eta$ on $\g$ taken from $\os$. This then obviously avoids the issue outlined above, as $\nabla^3\eta\in L^\infty(0,T;L^2(\os))$ from our framework introduced later in this Section.

Now, this does not avoid all difficulties. Once this formulation is adopted, a formidable difficulty arises when estimating terms such as
 $$\int_0^t\int_{\of}q \hd^2 \a_i^j  \hd^2 v,_j^i dx dt\,$$ since $q$ does not behave like a dissipative term. Although, as we will see, $\hd q$ is indeed dissipative, the fact $q$ is linked to $\eta,_3$ on $\g$ unfortunately does not allow any conclusion on making $q$ behave like $\hd q$.
 
 Instead, we will show the rather surprising property that the Stokes extension of $\hd\eta$ in $\of$ is actually a dissipative term, making the integral term above possible to estimate in a good way, independently of time.

 \subsection{An extension of $\eta$  (coming from the solid phase) into $\of$}

\subsubsection{Definition of our Arbitrary Lagrangian map}
We denote by $\e$ the (periodic in the horizontal directions) extension of $\eta$ (defined in $\os$) into $\of$ defined by the linear Stokes problem:
\begin{subequations}
	\label{he}
	\begin{alignat}{2}
		-\Delta\e+\nabla f&= 0 \ \  \ &&\text{in} \ \ \of \,,\\
		\div\e&=-\frac{1}{|\of|}\int_\g\eta^3 dx_h\ \  \ &&\text{in} \ \ \of \,,\label{he.b}\\
		\e &= \eta  \ \  &&\text{on} \ \ \g \,,\\
		\e &= 0\ \  && \text{on} \ \ \Gamma_{top} \,,
	\end{alignat}
\end{subequations}
which we then use in the formulation (\ref{ale}).

\begin{remark}
The compatibility condition (\ref{he.b}) ensures that the linear Stokes problem (\ref{he}) has a unique solution, with standard regularity properties, due to the necessary compatibility condition
$$\int_\of \div\e dx=\int_{\p\of} \e\cdot N^f dx_h\,.$$ 
\end{remark}

\begin{remark}
	The reason we use this Stokes extension and not an harmonic extension is for (\ref{2111.1}) established after, as it allows for $\hd \div\e=0$, whereas an harmonic extension would not provide this.
\end{remark}

We then see from the Dirichlet boundary condition on $\g$ in (\ref{he}) that with the Lagrangian displacement $\eta$ $$(\tilde\eta+\text{Id})(\R^2\times\{h_s,h\},t)=(\eta+\text{Id})(\R^2\times\{h_s,h\},t)\,,$$
and so $\tilde\eta+\text{Id}$ describes the motion of the fluid domain at each time.

\begin{remark}
Note that for the basic individual cell $\of$, although we may not have $(\tilde\eta+\text{Id})(\of,t)=(\eta+\text{Id})(\of,t)$, the thing that matters is that $(\tilde\eta+\text{Id})(\g,t)=(\eta+\text{Id})(\g,t)$.
\end{remark}

\subsubsection{Properties enabled by elliptic regularity for the linear extension operator}
\label{ellipticsub}

The first property needed later on is that for any first order horizontal derivative
\begin{align*}
	-\Delta(\hd\e)+\nabla(\hd f)&= 0 \ \  \ \text{in} \ \ \of \,,\\
		\div(\hd\e)&= 0 \ \  \ \text{in} \ \ \of \,,\\
	\hd\e  &= \hd\eta^s \ \  \text{on} \ \ \g \,,\\
	\hd\e &= 0\ \   \text{on} \ \ \Gamma_{top} \,,
\end{align*}
(we remind the top of the domain is horizontal for the Dirichlet condition on $\Gamma_{top}$), we have by standard elliptic regularity
\begin{equation}
	\label{te.2}
\|\hd\e\|_{H^2(\of)}\lesssim \|\hd\eta^s\|_{H^{\frac 3 2}(\g)}=\|\hd(\eta^s-\eta_e)\|_{H^{\frac 3 2}(\g)}\,,
\end{equation}
where we used $\hd\eta_e=0$.

Standard regularity for the first, second and third time-differentiated versions of (\ref{he}) yield

\begin{equation}
\label{tildev}
\|\v\|^2_{H^3(\of)}\lesssim \|v^f\|^2_{H^{\frac 5 2}(\g)}\lesssim D(t)^2 \,,
\end{equation}
\begin{equation}
	\label{tildevt}
	\|\v_t\|^2_{H^2(\of)}\lesssim \|v^f_{t}\|^2_{H^{\frac 3 2}(\g)}\lesssim D(t)^2 \,.
\end{equation}
\begin{equation}
	\label{tildevtt}
	\|\v_{tt}\|^2_{H^1(\of)}\lesssim \|v^f_{tt}\|^2_{H^{\frac 1 2}(\g)}\lesssim D(t)^2 \,.
\end{equation}

 where the dissipative norm $D(t)$ is defined in (\ref{D}). Note that for (\ref{tildev}), (\ref{tildevt}) and (\ref{tildevtt}), we take the viewpoint the trace is taken from the fluid phase $\of$, whereas for (\ref{te.2}) we take the trace from the solid phase $\os$.
 
  This is why $\tilde v$, $\tilde v_t$ and $\tilde v_{tt}$ have their norms controlled in $L^2(0,T; H^3(\of))$, $L^2(0,T; H^2(\of))$ and $L^2(0,T; H^1(\of))$ respectively, whereas $\tilde\eta$ does not seem to have this kind of control. As we will see later on, $\hd\tilde\eta$ is also controlled in  $L^2(0,T; H^2(\of))$ independently of $T$, for initial data close enough to equilibrium.
  
  We will also need the estimate for $\e$ when no horizontal derivative is applied.
  We first notice that elliptic regularity applied to (\ref{he}) yields
  \begin{equation*}
  \|\e\|_{H^3(\of)}\lesssim \left|\int_{\g} \eta^3 dx_h\right|+\|\eta\|_{H^{\frac 5 2}(\g)}\lesssim \|\eta-\eta_e\|_{H^{\frac 5 2}(\g)}\,,
  \end{equation*}
  where we used $\eta_e=0$ on $\g$ in the last inequality. This then implies for the norm $N$ defined in (\ref{N}):
  \begin{equation}
  	\label{te}
  	\|\e\|_{H^3(\of)}\lesssim N(t)\,.
  \end{equation}

\subsection{Local in time existence}

We now define the norm $N(t)$ which is finite initially, and remains finite (so long as the solution exists):

\begin{align}
	N(t)=& \sum_{k=0}^2 \|\partial_t^k v(\cdot,t)\|_{H^{2-k}(\Omega_0^f)}+ \sum_{k=0}^2 \|\partial_t^k v(\cdot,t)\|_{H^{2-k}(\Omega_0^s)}+\|\eta(\cdot,t)-\eta_e\|_{H^3(\os)}\,,
\label{N}
\end{align}
where $\eta_e$ is defined in (\ref{0911.1}), 
as well as the dissipative norm $D$ (associated to the viscous phase) which is square integrable in time (so long as the solution exists):
\begin{align}
	D(t)=& \sum_{k=0}^2 \|\partial_t^k v(\cdot,t)\|_{H^{3-k}(\of)}\,.
	\label{D}
\end{align}

This dissipative norm has a particularily important part to play in establishing global in time existence.



The total energy $E(t)$ is then defined as
\begin{equation}
E(t)=\sup_{[0,t]}N(t)^2+\int_0^t D(s)^2 ds\,.
\label{E}
\end{equation}

Due to the good elliptic properties of subsection \ref{ellipticsub}, the local in time existence of (\ref{ale}) follows directly from the known local in time existence of the problem in Lagrangian variables.

\begin{remark}
	Of course, $v_t(\cdot,0)$ and $v_{tt}(\cdot,0)$ are functions of the initial data $\eta(\cdot,0)$ and $v(\cdot,0)$, which means that $N(0)<\infty$ implies that $v(\cdot,0)$ is smoother than $H^2(\of)$ or $H^2(\os)$. See for instance the property (\ref{higher}) which establishes later that in fact $v\in L^\infty(0,T;H^3(\of))$ (although this property is not really needed in this paper). What matters is that the finite norm $N(t)$ is the framework giving existence of a solution to the problem which keeps this property propagating in (positive) time.
\end{remark}

\subsection{Small data assumption}

Therefore, given an initial data $(v_0(x), \eta_0(x))$ satisfying the appropriate compatibility conditions and regularity, we have that the problem (\ref{ale}) has a local in time solution on $[0,T]$, $(v,\eta)$ such that 
\begin{equation}
E(T)<\infty\,.
\end{equation}

In this paper, we are interested in solutions close enough to equilibrium initially, namely such that $$E(0)<\epsilon\,,$$ for $\epsilon>0$ small. We know from the local in time theory that for a large time $T(\epsilon)>1$ the solution will still exist at $T(\epsilon)$ and that $E(T(\epsilon))\lesssim\epsilon$.

Due to (\ref{te}), this implies that for our Stokes extension of $\eta$ into $\of$, 
\begin{equation}
	\forall t\in [0,T(\epsilon)]\,,\ \tilde J\in [\frac{99}{100},\frac{101}{100}]\,,\ \|\A-\text{Id}\|_{L^\infty(\of)}\le \frac{1}{100} \,.
	\label{tildeapriori}
\end{equation}
for $\epsilon>0$ small enough.

It is important to keep in mind these estimates, as they allow for a faster treatment of the global in time proof.


We will prove that for $\epsilon>0$ small enough, the solution is defined on $[0,\infty)$, its energy $E(t)$ remains small, and the solution converges in large time towards a flat interface solution, as defined in \ref{flat}.
  
 \subsection{Statement of the Theorems}\label{sec:maintheorem}
 
 \begin{theorem}[Global in time existence]\label{theorem_main} Let us assume that
 	\begin{enumerate}
 		\item  Our initial data satisfies the required compatibility conditions for local existence of a smooth solution (with finite $N(t)$).
 		
 		\item The elastic coefficient is large enough relative to gravity:
 		\begin{equation}\label{lambdalarge}
 			\d\lambda\ge \max(1, c g^2)
 		\end{equation} for $c>0$ large enough.
 		\item  The initial energy is small enough: $E(0)=N(0)^2\le\epsilon_0$\,,
 	\end{enumerate}
 	for $\epsilon_0>0$ small enough. Then, the local in time solution exists for all positive time and $E(t)$ remains small, of order $\epsilon_0$.
 	\end{theorem}

 	\begin{remark}
 	In the absence of gravity, $g=0$, and the condition (\ref{lambdalarge}) in Theorem \ref{theorem_main} reduces to $\lambda\ge 1$. This condition can be relaxed for $g\ne 0$ to $\lambda>0$, but for the sake of simplicity in the proofs, it is assumed this way.
 	\end{remark}
 	\begin{remark}
 		Due to the last term in the energy (\ref{N}), the smallness condition on $E(0)$ requires $h_e$ (average height of the elastic phase) to be close to $h_s$ (height of the reference domain for the elastic phase, where the wave equation is set).
 	\end{remark}
 	
 	\begin{theorem}[Asymptotic convergence]\label{theorem_asymptotic} Let us assume that
 	the assumptions of Theorem \ref{theorem_main} are satisfied.  Then,  the interface $\Gamma(t)$ converges towards the flat interface $\Gamma_e$ (defined in (\ref{equilibriumboundary})) in $H^{\frac 5 2}(\g)$, while the velocity in the fluid phase converges to zero as $t\rightarrow\infty$, and the vertical component of the displacement converges to the solution of a one dimensional wave equation in the solid phase. To be more specific, the norms in which these convergences hold are:
 	\begin{enumerate}
 		\item $\displaystyle\lim_{t\rightarrow\infty} \|v\|_{H^2(\of)}=0=\displaystyle\lim_{t\rightarrow\infty} \|v_t\|_{H^1(\of)}\,,$ 
 		\item   $\displaystyle\lim_{t\rightarrow\infty} \|\eta-(0, 0, h_e-h_s)\|_{H^{\frac{5}{2}}(\g)}=0\,,$
 		\item $ \displaystyle\lim_{t\rightarrow\infty} \|v^h\|_{H^1(\os)}=0=\displaystyle\lim_{t\rightarrow\infty} \|\eta^h\|_{H^2(\os)}$, where $f^h=(f_1,f_2)\,,$
 		\item There exists $(\alpha_0,\alpha_1)\in H^1_0(0,h_s)\times H^1_0(0,h_s)$ such that $\alpha$ defined by (\ref{wave1d}) satisfies
 		$$ \displaystyle\lim_{t\rightarrow\infty} \|v^3-\alpha_t\|_{L^2(\os)}=0=\displaystyle\lim_{t\rightarrow\infty} \|\eta^3-\alpha\|_{H^1(\os)}\,.$$
 	\end{enumerate}
 	
 	\begin{remark}
 		The existence of \it{flat interface solutions} as defined in \ref{flat} show that we can only hope to obtain the type of convergence as above, and that convergence towards the canonical equilibrium would not be possible for any solution satisfying Theorem \ref{theorem_main} (even in the case $h_e=h_s$).
 		
 	\end{remark}
 	
 	\end{theorem}

\subsection{Variational formulation of the problem and its first and formal second time differentiated versions in Arbitrary Lagrangian representation}

The variational formulation of the problem in these variables then becomes for all $t\in [0,T]$ and any $\phi\in H^1_0(\Omega;\R^3)$:
\begin{align}
	0=&\int_{\of} \J v_t\cdot \phi\ dx+ \int_{\of} (v-\tilde v)_l \a_l^j v,_j \cdot \phi\ dx +\nu\int_{\of} \a_j^l \A_j^k v,_k\cdot\phi,_l\ dx  -\int_{\of} \a_i^j q \phi^i,_j dx\n\\
	& + \int_{\os}v_t\cdot\phi\ dx+\lambda\int_{\os}\nabla \eta\cdot\nabla\phi\ dx+ g \int_\of \J \phi^3 dx+  g \int_\os \phi^3 dx
	\,.
	\label{var1}
\end{align}	

We will also need the first time differentiated version of this problem which tell us that 
 for all $t\in [0,T]$ and any $\phi\in H^1_0(\Omega;\R^3)$:
\begin{align}
	0=&\int_{\of} (\J v_{t})_t\cdot \phi\ dx+ \int_{\of} ((v-\tilde v)_l \a_l^j v,_j)_t \cdot \phi\ dx +\nu\int_{\of} (\a_j^l \A_j^k v,_k)_t\cdot\phi,_l\ dx -\int_{\of} (\a_i^j q)_t \phi^i,_j dx\n\\
	 &+ \int_{\os}v_{tt}\cdot\phi\ dx+\lambda\int_{\os}\nabla v\cdot\nabla\phi\ dx+ g \int_\of \J_t \phi^3 dx
\,,
	\label{var2}
\end{align}	

For the second time-differentiated version of this problem, the energy inequality satisfied for $v_{tt}$ at this level is obtained as if we formally had the formal time-differentiated problem of (\ref{var2}) satisfied (for $v_{tt}$ as test function):
\begin{align}
	0=&\int_{\of} (\J v_{t})_{tt}\cdot \phi\ dx+ \int_{\of} ((v-\tilde v)_l \a_l^j v,_j)_{tt} \cdot \phi\ dx +\nu\int_{\of} (\a_j^l \A_j^k v,_k)_{tt}\cdot\phi,_l\ dx
	 -\int_{\of} (\a_i^j q)_{tt} \phi^i,_j dx\n\\
	 & + \int_{\os}v_{ttt}\cdot\phi\ dx+\lambda\int_{\os}\nabla v_t\cdot\nabla\phi\ dx+ g \int_\of \J_{tt} \phi^3 dx
	\,.
	\label{var3}
\end{align}

As we will see later in (\ref{var1h3}), a term arising in the study of (\ref{var1}) is $$\int_0^t\int_{\of}q \hd^2 \a_i^j  \hd^2 v,_j^i dx dt\,.$$

The issue to treat this space time integral is that $q$ is linked to $\eta,_3^3$ on $\Gamma$, which is a priori not a good term for $L^2$ in time energies as needed in the integral, where a priori the only good term is the term in derivatives of $v$.

To get around this seemingly impossible to resolve issue, we will establish in the Section \ref{crucial} that $\hd\tilde\eta$ is controlled in $L^2(0,t;H^2(\of))$ independently of time by a small constant (in small data regime) multiplying $\sqrt{E(t)}$ (and can thus be viewed as a dissipative term, like $\nabla^3 v$ in $\of$), plus an initial term. This is not an easy result, as the functional framework provides us a priori just with an $L^\infty$ in time estimate for $\tilde\eta$, which would mean a priori  $\hd\tilde\eta$ would be controlled in $L^2(0,t;H^2(\of))$ with a $\sqrt{t}$ growth, unsuitable for large time.

	Before we can move to this crucial Section \ref{crucial}, we first establish hereafter a number of properties that will be needed later on in the paper for $q$ and $q_t$. Although only (\ref{nablaq}) and (\ref{0511.2}) are needed for Section \ref{crucial}, it is more convenient to have all needed relations later on for the pressure in the same Section \ref{sectionpressure}.

\section{Pressure estimates in terms of $D$ and $N$}\label{sectionpressure}

\subsection{High order estimates for $q$ and $q_t$ in $\of$}

Our starting point is the identity $$\nabla\tilde X \A=\text{Id}\,,$$ which implies
$$\tilde X,_k^i \A_i^m=\delta_k^m\,,$$

Multiplying the Arbitrary Lagrangian Navier-Stokes equation (\ref{ale.a}) by $\frac{\tilde X,_m^i}{\J}=\frac{\tilde\eta,_m^i+\delta_m^i}{\J}$ (for $m$ fixed), and summing over all $i$, yields:
\begin{equation}
	\label{qale}
	(v_t^i+(v-\tilde v)_l \A_l^j v,_j^i-\nu \A_j^l (\A_j^k v^i,_k),_l+g\delta_3^i)(\e,_m^i+\delta_m^i)=-q,_m\,.
\end{equation}

Applying (\ref{qale}) to $m=\alpha\in\{1,2\}$, we have for the order $1$ gravity term $\delta_3^i \delta_\alpha^i=0$, which implies using $\|ab\|_{H^1(\of)}\lesssim \|a\|_{H_2(\of)}\|b\|_{H^1(\of)}$ that
\begin{align*}
	\|\hd q\|_{H^1(\of)}\lesssim &   \|v_t\|_{H^1(\of)}\|\nabla\e+\text{Id}\|_{H^2(\of)}+  \|v-\tilde v\|_{H^2(\of)}\|\nabla v\|_{H^1(\of)}\|\A\|_{H^2(\of)}\n\\
	&+   \|\A\|_{H^2(\of)}\|\A_j^k v^i,_k\|_{H^2(\of)}+ g \|\hd\e^3\|_{H^1(\of)}\n\\
	\lesssim &   \|v_t\|_{H^1(\of)}(\|\e\|_{H^3(\of)}+1) +  \|v-\tilde v\|_{H^2(\of)}\|\nabla v\|_{H^1(\of)}\|\A\|_{H^2(\of)}\n\\
	&+   \|\A\|^2_{H^2(\of)}\|\nabla v\|_{H^2(\of)}+ g \|\hd\e^3\|_{H^1(\of)}
	\,,
\end{align*}	
where we used $H^2(\of)$ is a Banach algebra in the second to last term of the second inequality above.

Using our close to identity assumption (\ref{tildeapriori}), we see that
\begin{equation*}
	\|\nabla^2 \A\|_{L^2(\of)}\lesssim \|\nabla^3\e\|_{L^2(\of)}+\|\nabla^2\e\|^2_{L^4(\of)}\lesssim \|\e\|_{H^3(\of)}+\|\e\|^2_{H^3(\of)}\,,
\end{equation*}
where we used Sobolev's embeddings in the second inequality above. Using this for the previous inequality for $\hd q$ yields

\begin{align*}
	\|\hd q\|_{H^1(\of)}\lesssim &   \|v_t\|_{H^1(\of)}(\|\e\|_{H^3(\of)}+1)+  \|v-\tilde v\|_{H^2(\of)}\|\nabla v\|_{H^1(\of)}(1+\|\e\|^2_{H^3(\of)})\n\\
	&+   (1+\|\e\|^4_{H^3(\of)})\|\nabla v\|_{H^2(\of)}+ g \|\hd\e^3\|_{H^1(\of)}\,.
\end{align*}	
Using (\ref{te}) and (\ref{tildev}), we then infer
\begin{equation}
	\|\hd q\|_{H^1(\of)}\lesssim   D(t) (1+P(N(t))  ) + g\|\hd\e^3\|_{H^1(\of)}
	\,,\label{nablaq}
\end{equation}
where we remind our convention on $P$ is in subsection \ref{polynomials}.

The same application of (\ref{qale}) yields similarly
\begin{equation}
	\|\hd q\|_{H^{\ud}(\g)}\lesssim D(t) (1+P(N(t))) + g\|\hd\e^3\|_{H^{\ud}(\g)}\,.\label{0511.2}
\end{equation}

Applying  (\ref{qale}) to $m=3$ yields similarly
\begin{equation}
	\|q,_3+g\|_{H^1(\of)} \lesssim  D(t) (1+P(N(t))) + g\|\e,_3^3\|_{H^1(\of)}\,.\label{1307.1}
\end{equation}
Due to (\ref{te}), we then infer from (\ref{1307.1}) that
\begin{equation}
	\|q,_3+g\|_{H^1(\of)} \lesssim  D(t) (1+P(N(t))) + g N(t)\,,\label{1211.1}
\end{equation}
which with (\ref{nablaq}) and definition (\ref{0911.2}) allows
\begin{equation}
	\|\nabla (q-q_e)\|_{H^1(\of)}\lesssim D(t) (1+P(N(t))  ) +  g N(t)\,.\label{1307.3}
\end{equation}

\begin{remark} For the case without gravity, this relation shows the gradient of $q$ will be square integrable in time. We will also show the second term on the right hand side of (\ref{nablaq}) to be square integrable in time, if $\lambda>0$ is large enough relative to $g$.
\end{remark}

Similarly, taking one time derivative of (\ref{qale}) (which has for effect to replace the linear term in $\e$ by a linear term in $\v$, which is dissipative due to (\ref{tildev})) also yields:
\begin{align}
	\|\nabla q_t\|_{L^2(\of)}\lesssim    D(t) (1+P(N(t)))
	\,.\label{nablaqt}
\end{align}



By Poincar\'e-Wirtinger, our control (\ref{nablaqt}) implies
\begin{equation}
	\|q_t-\frac{1}{|\g|}\int_{\g} q_t\ dx_h\|_{L^2(\of)}\lesssim   D(t) (1+P(N(t)))
	\,.\label{qt}
\end{equation}

In a way similar as the standard Poincar\'e-Wirtinger inequality is proved in a compact domain, taking a sequence $f_n$ in $H^2(\of)$ such that $\|f_n\|_{H^2(\of)}=1$ with $\d\int_\g f_n dx_h=0\,,$ and $\d\|\nabla f_n\|_{H^1(\of)}\le \frac{1}{n} \|f_n\|_{L^\infty(\of)}\,,$ leads to a contradiction. This therefore establishes that there exists $C>0$ finite such that
\begin{equation}
	\forall f\in H^2(\of)\,,\ \|f-\frac{1}{|\g|}\int_{\g} f\ dx_h\|_{L^\infty(\of)}\le C \|\nabla f\|_{H^1(\of)}\,.
	\label{newpw}
\end{equation}
Using (\ref{newpw}) for $q-q_e$ and (\ref{1307.3}) we then infer
\begin{equation}
	\|q-q_e-\frac{1}{|\g|}\int_{\g} q-q_e\ dx_h\|_{L^\infty(\of)}\lesssim   D(t) (1+P(N(t)))+ g N(t)
	\,.\label{qinfty}
\end{equation}

\subsection{Estimates on $\g$ for $q$ and $q_t$}

Multiplying the continuity of stress (\ref{systemhorizontal}) by $\tilde X,_3^i=\tilde \eta,_3^i+\delta_3^i$, summing over all $i$, and using $N=(0,0,-1)$ on $\g$ yields:
\begin{equation}
	q\ \text{det}\nabla \tilde X=(-\lambda\eta,_3^i+\nu \a_j^3 \A_j^k v^i,_k)(\e,_3^i+\delta_3^i)  \,,
	\label{qbc}
\end{equation}
which then provides (so long as (\ref{tildeapriori}) is satisfied):
\begin{equation}
	\|q\|_{L^4(\g)}\le  C_\lambda P(N(t))=P_\lambda(N(t))\,.\label{177.2bis}
\end{equation}

Taking one time derivative of (\ref{qbc}) then shows that $q_t$ on $\g$ is a sum of product terms where appears either $\nabla v^f$, $\nabla v^s$, $\nabla v^f_t$ or $\nabla \tilde v$, which yields:
\begin{equation}
\|q_t\|_{L^2(\g)}\le  C_\lambda(N(t)+D(t))(1+P(N(t)))\,.\label{177.2ter}
\end{equation}

As a corollary of (\ref{nablaqt}), (\ref{qt}), and (\ref{177.2ter}), we have by Sobolev embeddings that
\begin{equation}
\|q_t\|_{L^4(\of)}\le C_\lambda (D(t)+N(t))(1+P(N(t)))\,.\label{qtl4}
\end{equation}

\subsection{$L^\infty$ in time low order space estimates for $q$ and $q_t$ in $\of$}

With $\phi_0\in H^1_0(\Omega)$ such that $\div\phi_0= q$ in $\of$ and $\|\phi_0\|_{H^1_0(\Omega)}\lesssim \|q\|_{L^2(\of)}$ (see for instance Lemma 13 of \cite{CS1} for instance for a proof of the existence of such $\phi_0$), we obtain by using $\phi_0$ in (\ref{var1}) that
\begin{align*}
\|q\|^2_{L^2(\of)}- \|q\|^2_{L^2(\of)}\|\tilde a-\text{Id}\|_{L^\infty(\of)}\lesssim&  \frac{1}{\epsilon} P_\lambda(N(t))^2+\frac{1}{\varepsilon} g^2 (1+\|\tilde J\|^2_{L^2(\of)})+ \varepsilon \|\phi_0\|^2_{H^1(\of)}\n\\
\lesssim& \frac{1}{\epsilon} P_\lambda(N(t))^2+\frac{1}{\varepsilon} g^2(1+\|\tilde J\|^2_{L^2(\of)}) + \varepsilon \|q\|^2_{L^2(\of)}
\,,
\end{align*}
which with our assumption (\ref{tildeapriori}) implies (for $\varepsilon>0$ chosen small enough)
\begin{equation*}
\|q\|^2_{L^2(\of)}\lesssim  P_\lambda(N(t))^2+g^2\,.
\end{equation*}

With $\phi_1\in H^1_0(\Omega)$ such that $\div\phi_1= q_t$ in $\of$ and $\|\phi_1\|_{H^1_0(\Omega)}\lesssim \|q_t\|_{L^2(\of)}$, we obtain by using $\phi_1$ in (\ref{var2}) that
\begin{align*}
\|q_t\|^2_{L^2(\of)}- \|q_t\|^2_{L^2(\of)}\|\tilde a-\text{Id}\|_{L^\infty(\of)}\lesssim&    \frac{1}{\epsilon} P_\lambda(N(t))^2+\frac{1}{\varepsilon} g^2\|\tilde J_t\|^2_{L^2(\of)}+ \varepsilon \|\phi_1\|^2_{H^1(\of)}\n\\
\lesssim& \frac{1}{\epsilon} P_\lambda(N(t))^2+\frac{1}{\varepsilon} g^2 \|\tilde J_t\|^2_{L^2(\of)} + \varepsilon \|q_t\|^2_{L^2(\of)}
\,,
\end{align*}
which with our assumption (\ref{tildeapriori}) implies (for $\varepsilon>0$ chosen small enough)
\begin{equation}
\|q_t\|^2_{L^2(\of)}\lesssim  P_\lambda(N(t))^2\,.\label{qmultiplier}
\end{equation}

\begin{remark}
We have higher regularity properties for $q$ and $q_t$, see for instance (\ref{higher}) later. These properties are however not needed in this paper.
\end{remark}

\section{A crucial $L^2$ in time estimate for $\hd^2\eta$ on $\g$}\label{crucial}

We establish in this Section the fundamental inequality (\ref{inter15}), from which follows the needed estimate:
\begin{equation}
	\int_0^t\int_\of |\hd^2 \a_i^j|^2 dx dt\lesssim (1+\sup_{[0,t]} P(N)) E(t) + C_\lambda E(0)\,.\label{177.1}
\end{equation} 


\subsection{$L^2$ in time estimate on $\g$ for the wave equation}

Our starting point is to take horizontal derivatives of the linear wave equation, take the scalar product with two horizontal derivatives and one vertical derivative of $\eta$ and integrate in space-time. As we seek to obtain fractional derivative space regularity on $\g$, it is actually very convenient to work with Fourier series as in our reminder in Section \ref{sec:sobolev-spaces-gamma}.

From our reminder in Section \ref{sec:sobolev-spaces-gamma}, since by horizontal periodicity $$\int_{[0,L]^2} \hd\eta(x_h,x_3,t) dx_h=0\,,$$
we have expanding $\hd\eta(\cdot,x_3,t)$ (where $\hd$ denote here any of the first order horizontal derivatives with respect to $x_1$ or $x_2$ ) in this basis that for each component $i=1, 2, 3$,
\begin{equation}
	\hd\eta^i(x_h,x_3,t)=\sum_{n=1}^\infty x^n(x_3,t) e_n(x_h)\,,\label{070824.1}
\end{equation}
with $$x^n(x_3,t)=\int_{[0,L]^2} \hd\eta^i(x_h,x_3,t) e_n(x_h) dx_h\,.$$

Since $\hd\eta$ is solution of the wave equation (without gravity term), we have for each mode (using $\Delta_0 e^n=-\lambda_n e^n$):
\begin{equation}
	x^n_{tt}+\lambda\lambda_n x^n-\lambda x,_{33}^n=0\,.\label{287.1}
\end{equation}
Multiplying (\ref{287.1}) by $x,_3^n$, and integrating over $[0,h_s]\times [0,t]$ yields:
\begin{align*}
	0=&\int_0^t\int_0^{h_s} x^n_{tt} x,_3^n +\lambda\lambda_n x^n x,_3^n -\lambda x,_{33}^n x,_3^n dx_3 dt\n\\
	=&\int_0^t\int_0^{h_s} - x^n_{t} x_t,_3^n +\lambda\lambda_n x^n x,_3^n -\lambda x,_{33}^n x,_3^n dx_3 dt+\left[\int_0^{h_s} x^n_t x,_3^n dx_3\right]_0^t\n\\
	=&\ud\int_0^t\int_0^{h_s} - |x^n_{t}|^2,_3 +\lambda\lambda_n |x^n|^2,_3 -\lambda|x,_{3}^n|^2,_3 dx_3 dt+\left[\int_0^{h_s} x^n_t x,_3^n dx_3\right]_0^t\,.
\end{align*}
Integrating by parts in space then yields
\begin{align}
	0
	=&\ud\int_0^t (- |x^n_{t}|^2 +\lambda\lambda_n |x^n|^2 -\lambda|x,_{3}^n|^2)(h_s,\cdot)  dt-\ud\int_0^t (- |x^n_{t}|^2 +\lambda\lambda_n |x^n|^2 -\lambda|x,_{3}^n|^2)(0,\cdot)  dt\n\\
	&+\left[\int_0^{h_s} x^n_t x,_3^n dx_3\right]_0^t\,.\label{287.2}
\end{align}
From the Dirichlet boundary condition on $\Gamma_B$, $x^n=0=x^n_t$ on $x_3=0$, and so (\ref{287.2}) implies:
\begin{align*}
	\int_0^t \lambda\lambda_n |x^n|^2 (h_s,\cdot) dt
	=&\int_0^t ( |x^n_{t}|^2 +\lambda|x,_{3}^n|^2)(h_s,\cdot)  dt-\int_0^t \lambda|x,_{3}^n|^2(0,\cdot) dt-2\left[\int_0^{h_s} x^n_t x,_3^n dx_3\right]_0^t\n\\
	\le& \int_0^t  (|x^n_{t}|^2 +\lambda|x,_{3}^n|^2)(h_s,\cdot) dt-2\left[\int_0^{h_s} x^n_t x,_3^n dx_3\right]_0^t\,.
\end{align*}
We multiply this identity by $\lambda_n^{\ud}$ and sum over $n$:
\begin{align*}
	\sum_{n=1}^\infty \int_0^t\lambda \lambda_n^{\frac 3 2} |x^n|^2(h_s,\cdot)   dt
	\le& \sum_{n=1}^\infty \int_0^t (\lambda_n^{\frac 1 2}|x^n_{t}|^2 +\lambda\lambda_n^{\frac 1 2}|x,_{3}^n|^2)(h_s,\cdot) dt-2\sum_{n=1}^\infty \left[\int_0^{h_s}\lambda_n^{\frac 1 2} x^n_t x,_3^n dx_3\right]_0^t\,.
\end{align*}
Therefore,
\begin{align}
	\lambda\int_0^t \|\hd \eta^i\|_{H^{\frac 3 2}(\g)}^2 dt
	\le& \int_0^t  \|\hd v^i\|_{H^{\frac 1 2}(\g)}^2 +  \lambda\|\hd \eta^i,_3\|_{H^{\frac 1 2}(\g)}^2 dt-2  \sum_{n=1}^\infty \left[\int_0^{h_s}\lambda_n^{\frac 1 2} x^n_t x,_3^n dx_3\right]_0^t\n\\
	\le& \int_0^t  \|\hd v^i\|_{H^{\frac 1 2}(\g)}^2 + \lambda \|\hd \eta^i,_3\|_{H^{\frac 1 2}(\g)}^2 dt+   \sum_{n=1}^\infty \lambda_n^{\frac 1 2} \int_0^{h_s} |x^n_t|^2(t)+| x,_3^n|^2(t) dx_3 \n\\
	&+ \sum_{n=1}^\infty \lambda_n^{\frac 1 2} \int_0^{h_s} |x^n_t|^2(0)+| x,_3^n|^2(0) dx
	\,.\label{287.3}
\end{align}
Next we see by multiplying (\ref{287.1}) by $x^n_t$ and integrating in $[0,h_s]\times [0,t]$ that
\begin{align*}
	0=&\int_0^t\int_0^{h_s} x^n_{tt} x^n_t +\lambda\lambda_n x^n x^n_t -\lambda x,_{33}^n x^n_t dx_3 dt\n\\
	=&\int_0^t\int_0^{h_s} x^n_{tt} x^n_t +\lambda\lambda_n x^n x^n_t +\lambda x,_{3}^n {x,_3^n}_t dx_3 dt-\lambda \int_0^t\left[x,_3^n x^n_t\right]_0^{h_s}  dt\n\\
	=&\ud \left[\int_0^{h_s} | x^n_t|^2 +\lambda\lambda_n |x^n|^2 +\lambda|x,_{3}^n|^2 dx_3\right]_0^t-\lambda\int_0^t x,_3^n x^n_t (h_s,\cdot)  dt\,.
\end{align*}
Multiplying this identity by $\lambda_n^{\frac 1 2}$, and summing over $n$ yields:
\begin{align*}
	\sum_{n=1}^\infty \lambda_n^{\frac 1 2} \int_0^{h_s} |x^n_t|^2(t)+\lambda| x,_3^n|^2(t) dx_3=& -\sum_{n=1}^\infty \lambda_n^{\frac 1 2} \int_0^{h_s} \lambda\lambda_n |x^n|^2(t) dx_3\n\\
	&+  \sum_{n=1}^\infty \lambda_n^{\frac 1 2} \int_0^{h_s} |x^n_t|^2(0)+\lambda\lambda_n |x^n|^2(0) +\lambda| x,_3^n|^2(0) dx_3\n\\
	&+2\sum_{n=1}^\infty \lambda_n^{\frac 1 2}\lambda \int_0^t x,_3^n x^n_t (h_s,\cdot) dt\n\\
	\le & \sum_{n=1}^\infty \lambda_n^{\frac 1 2} \int_0^{h_s} |x^n_t|^2(0)+\lambda\lambda_n |x^n|^2(0) +\lambda| x,_3^n|^2(0) dx_3\n\\
	&+\sum_{n=1}^\infty \lambda_n^{\frac 1 2} \int_0^t (\lambda^2|x,_3^n|^2+| x^n_t|^2) (h_s,\cdot)dt\n\\
	\le & \sum_{n=1}^\infty \lambda_n^{\frac 1 2} \int_0^{h_s} |x^n_t|^2(0)+\lambda\lambda_n |x^n|^2(0) +\lambda| x,_3^n|^2(0) dx_3\n\\
	&+\int_0^t\lambda^2\|\hd\eta,_3\|_{H^{\frac 1 2}(\g)}^2+\|\hd v\|_{H^{\frac 1 2}(\g)}^2  dt\,.
\end{align*}
Since $\lambda_n\ge \lambda_1\ge 1$, we infer from this inequality that
\begin{align}
	\sum_{n=1}^\infty \lambda_n^{\frac 1 2} \int_0^{h_s} |x^n_t|^2(t)+\lambda| x,_3^n|^2(t) dx_3
	\le & C \sum_{n=1}^\infty \lambda_n \int_0^{h_s} |x^n_t|^2(0)+\lambda\lambda_n |x^n|^2(0) +\lambda| x,_3^n|^2(0) dx_3\n\\
	&+\int_0^t \lambda^2\|\hd\eta,_3^i\|_{H^{\frac 1 2}(\g)}^2+\|\hd v^i\|_{H^{\frac 1 2}(\g)}^2  dt\n\\
	\le & C \int_0^{h_s} \|\hd v^i\|_{H^1([0,L]^2\times\{x_3\})}^2(0)+\lambda \|\hd\eta^i\|_{H^2([0,L]^2\times\{x_3\})}^2(0) dx_3\n\\
	&+C\lambda\int_0^{h_s}\| \hd \eta,_3^i\|_{H^1([0,L]^2\times\{x_3\})}^2(0) dx_3
	+\int_0^t \lambda^2\|\hd\eta,_3^i\|_{H^{\frac 1 2}(\g)}^2+\|\hd v^i\|_{H^{\frac 1 2}(\g)}^2  dt\n\\
	\le &{(C+C\lambda)} E(0)+\int_0^t \lambda^2\|\hd\eta,_3^i\|_{H^{\frac 1 2}(\g)}^2+\|\hd v^i\|_{H^{\frac 1 2}(\g)}^2  dt\label{287.4}
	\,.
\end{align}

Using (\ref{287.4}) in (\ref{287.3}), and $\lambda\ge 1$, we obtain for each $i\in\{1,2,3\}$:
\begin{align}
	\int_0^t \lambda\|\hd \eta^i\|_{H^{\frac 3 2}(\g)}^2 dt\le&  C_\lambda E(0)+\int_0^t \lambda\|\hd\eta,_3^i\|_{H^{\frac 1 2}(\g)}^2+\|\hd v^i\|_{H^{\frac 1 2}(\g)}^2  dt+\int_0^t\lambda^2 \|\hd\eta,_3^i\|_{H^{\frac 1 2}(\g)}^2+\|\hd v^i\|_{H^{\frac 1 2}(\g)}^2  dt\n\\
	 \le&  C_\lambda E(0)+2\int_0^t (\lambda^2+\lambda)\|\hd\eta,_3^i\|_{H^{\frac 1 2}(\g)}^2+\|\hd v^i\|_{H^{\frac 1 2}(\g)}^2  dt\,.\label{inter6}
\end{align}

\subsection{Crucial $L^2$ in time estimate for $\|\hd\tilde\eta\|^2_{H^{2}(\of)}$ and  $\|\hd q\|^2_{H^{1}(\of)}$ :}

We remind that $\hd$ stands for any first order horizontal derivative with respect to $x_1$ or $x_2$. We begin by summing the identities (\ref{inter6}) from $1$ to $2$ (while omitting the summation symbols for conciseness), which means that hereafter the symbol $\hd$ is viewed as the horizontal gradient.

Using continuity of normal stress alongside $\Gamma$ (\ref{systemhorizontal}), (\ref{inter6}) implies for each $i\in\{1,2,3\}$:
\begin{align*}
	\int_0^t   \lambda\|\hd \e^i\|_{H^{\frac 3 2}(\g)}^2 dt \le &
	\int_0^t C \| v^i\|^2_{H^2(\of)}+ {2}(1+\frac{1}{\lambda})\|\hd(q \a_i^3 -\nu \a_l^3 \A_l^k v^i,_k) \|^2_{H^{\frac 1 2}(\g)}  dt
	+ C_\lambda E(0)\n\\
	\le &  \int_0^t  C\| v\|^2_{H^2(\of)}+ 4\|\hd q \a_i^3 +q \hd \a_i^3- \nu\hd(\a_l^3 \A_l^k v^i,_k) \|^2_{H^{\frac 1 2}(\g)}  dt
	+ C_\lambda E(0)\,,
\end{align*}	
as we assumed $\lambda\ge 1$.
Thus,

\begin{align}
\int_0^t   \lambda\|\hd \e^i\|_{H^{\frac 3 2}(\g)}^2 dt	
	\le &  \int_0^t  C\| v\|^2_{H^2(\of)}+  16\|\hd q \a_i^3\|^2_{H^{\frac 1 2}(\g)} + 16\nu^2\|\hd(\a_l^3 \A_l^k v^i,_k) \|^2_{H^{\frac 1 2}(\g)}+8\|q \hd \a_i^3 \|^2_{H^{\frac 1 2}(\g)}   dt\n\\
	&+ C_\lambda E(0)
	\,.\label{inter7}
\end{align}
Next, using (\ref{0211.1}) (vertical component) and (\ref{1407.1}), we have on $\g$
\begin{align*}
	\|q \hd \a_i^3 \|_{H^{\frac 1 2}(\g)}\le& \|(q\underbrace{+\frac{gh_s}{2}}_{-q_e}) \hd \a_i^3 \|_{H^{\frac 1 2}(\g)}+\frac{gh_s}{2} \|\hd\a_i^3 \|_{H^{\frac 1 2}(\g)}\n\\
	\le& (C \|q-q_e\|_{H^{\frac 3 2}(\g)}+\frac{gh_s}{2})\|\hd\a_i^3 \|_{H^{\frac 1 2}(\g)}\n\\
	\le&  (C (\|v\|_{H^3(\of)}+\|\hd\tilde\eta\|_{H^{2}(\of)}+\lambda\|\eta-\eta_e\|_{H^3(\os)})(1+P(N)) + \frac{gh_s}{2} )\|\hd\a_i^3 \|_{H^{\frac 1 2}(\g)}\,,
\end{align*}
where we used for $q-q_e$ in (\ref{0211.1}) the fact that $$\a_3^3-\delta_3^3=(\tilde X,_1\times \tilde X,_2)^3-1\,,$$ which is a sum of first order horizontal derivatives and of products of such derivatives for $\e$. We also used (\ref{tildeapriori}).
Therefore,
	\begin{align}
		\|q \hd \a_i^3 \|_{H^{\frac 1 2}(\g)}
	\le& (C D(t) (1+P(N(t))) + P_\lambda(N(t))+ \frac{gh_s}{2}) \|\hd\a_i^3 \|_{H^{\frac 1 2}(\g)}\n\\
	\le&   ( P_\lambda(N(t)) +\frac{gh_s}{2}) \|\hd\e \|_{H^{\frac 3 2}(\g)} (1+P(N(t)))+ CD(t) P(N(t))\n\\
	\le& ( P_\lambda (N(t)) + \frac{gh_s}{2}) \|\hd\e \|_{H^{\frac 3 2}(\g)}+ CD(t) P(N(t))
	\,.\label{inter8}
\end{align}

Using (\ref{inter8}) and (\ref{0511.2}) for the $\hd q$ term in (\ref{inter7}) then yields (for $C>0$ independent of $t$, $\lambda$, and the initial data)
\begin{equation*}
	\int_0^t   \|\hd\tilde\eta\|^2_{H^{\frac 3 2}(\g)}  (\lambda- C{ g^2}-  P_\lambda (N(s))) ds \lesssim \int_0^t D^2(s) (1+P(N(s))) ds
	+ C_\lambda E(0)
	\,.
\end{equation*}
So long as $N\le 1$ on $[0,t]$, we then obtain the existence of some constant $C_1>0$ and $C_{0,\lambda}>0$ independent of $t$, and of the initial data, such that
\begin{equation*}
	\int_0^t   \|\hd\tilde\eta\|^2_{H^{\frac 3 2}(\g)}  (\lambda- C_1{ g^2}-  C_{0,\lambda} N(s)) ds \lesssim \int_0^t D^2(s) (1+P(N(s))) ds
	+ C_\lambda E(0)
	\,,
\end{equation*}
which thanks to our assumption (\ref{lambdalarge}) of a large enough $\lambda$ relative to $g^2$, and so long as our assumption (we can always assume $C_{0,\lambda}$ very large): 
\begin{equation}
	\label{smalldata}
N\le \frac{\lambda}{2 C_{0,\lambda}}\ll 1\,,
\end{equation}
 holds, implies
\begin{equation}
	\label{inter14}
	\int_0^t  \frac{\lambda}{4} \|\hd\tilde\eta\|^2_{H^{\frac 3 2}(\g)}   dt \lesssim  \int_0^t  D^2(s) (1+P(N(s))) ds 
	+ C_\lambda E(0)
	\,.
\end{equation}

\begin{remark}
We also assume to have chosen $C_{0,\lambda}$ large enough so that if (\ref{smalldata}) is satisfied, then our assumption (\ref{tildeapriori}) is also satisfied.
\end{remark}

By elliptic regularity (\ref{te.2}), this implies in turn
\begin{align}
\lambda	\int_0^t   \|\hd\tilde\eta\|^2_{H^{2}(\of)}   dt \lesssim&   \int_0^t D^2(s) (1+P(N(s)))  ds
	+ C_\lambda E(0)\n\\
	\lesssim&   (1+\sup_{[0,t]} P(N))\ \int_0^t D^2(s)  ds
	+ C_\lambda E(0)\n\\
	\lesssim&   (1+\sup_{[0,t]} P(N))\ E(t)
	+ C_\lambda E(0) \label{inter15}	
	\,.
\end{align}

This is  a crucial identity, as it allows us to state that the $L^2$ in time norm of $ \|\hd\tilde\eta\|_{H^{2}(\of)} $ is controlled by some $L^2$ in time norm of $v$ and $v_t$ in the fluid phase (for which dissipation provides natural control), and an extra constant from the initial data. Note also this result does not imply that $\|\hd\eta\|_{H^{2}(\os)} $ satisfies the same type of estimate.

Also note that in (\ref{inter15}), the constant multiplying $(1+\sup_{[0,t]} P(N))\ E(t)+C_\lambda E(0)$ in the $\lesssim$ convention does not depend on $\lambda$ (which is assumed large relative to $g$).

As a consequence of (\ref{inter15}), we infer from (\ref{nablaq}) that
\begin{equation}
	\int_0^t   \|\hd q\|^2_{H^{1}(\of)}   dt \lesssim  (1+\sup_{[0,t]} P(N))\ E(t)
	+ C_\lambda E(0)\,.\label{0611.2}
\end{equation}

\begin{remark}
	\label{tracksize}
	Note that (\ref{0611.2}) should have a $\lambda^{-1}$ factoring the right hand side, but since we assume $\lambda\ge 1$, we do not have to keep track of this term.
	
	This convention is not used in (\ref{inter15}), as we need to keep the $\lambda$ appearing there on the left hand side for some estimates.
	
\end{remark}

\section{Highest order in space estimate from the variational formulation}
\label{highspace}

We use $\hd^2\hd^2 v\in H^1_0(\Omega)$ (where $\hd^2$ denote any second order partial derivative $\frac{\partial^2}{\partial x_\alpha\partial x_\alpha}$, for $\alpha=1,2$) as test function in (\ref{var1}) and integrate by parts with respect to the horizontal variable twice to obtain:
\begin{align*}
	0=&\int_{\of} \hd^2 (\J v_t)\cdot \hd^2 v\ dx+ \int_{\of} \hd^2 ((v-\tilde v)_l \a_l^j v,_j) \cdot \hd^2 v\ dx +\nu\int_{\of} \hd^2(\a_j^l \A_j^k v,_k)\cdot\hd^2 v,_l\ dx\n\\
	& -\int_{\of} \hd^2(\a_i^j q) \hd^2 v,_j^i dx+ \int_{\os}\hd^2 v_t\cdot\hd^2 v\ dx+\lambda\int_{\os}\nabla \hd^2\eta\cdot\nabla\hd^2 v\ dx+ g \int_\of \hd^2\J \hd^2v^3 dx\n\\
	=&\int_{\of} \J \hd^2 v_t\cdot \hd^2 v\ dx+ \int_\of 2 \hd \J \hd v_t \cdot \hd^2 v + \hd^2\J v_t\cdot\hd^2 v dx+\int_{\of} \hd^2 ((v-\tilde v)_l \a_l^j v,_j) \cdot \hd^2 v\ dx\n\\
	& +\nu\int_{\of} \hd^2(\a_j^l \A_j^k v,_k)\cdot\hd^2 v,_l\ dx -\int_{\of} \hd^2(\a_i^j q) \hd^2 v,_j^i dx+ \int_{\os}\hd^2 v_t\cdot\hd^2 v\ dx+\lambda\int_{\os}\nabla \hd^2\eta\cdot\nabla\hd^2 v\ dx\n\\
	&+ g \int_\of \hd^2\J \hd^2v^3 dx
	\,.
\end{align*}

Integrating in time from $0$ to $t$ yields, and summing over $\alpha=1,2$ (while omitting to write the summation symbols for conciseness)  yields
\begin{align}
	&\ud\int_{\of}\J |\hd^2 v|^2(\cdot,t)\ dx \u{-\ud\int_0^t\int_\of \J_t |\hd^2 v|^2 dxdt}_{I_1}  +\u{\int_0^t\int_\of 2 \hd \J \hd v_t \cdot \hd^2 v + \hd^2\J v_t\cdot\hd^2 v dx dt}_{I_2}\n\\
	&+ \u{\int_0^t \int_{\of} \hd^2 ((v-\tilde v)_l \a_l^j v,_j) \cdot \hd^2 v\ dx dt}_{I_3}+ \nu\int_0^t \int_{\of} |\nabla\hd^2 v|^2 \ dx dt+ \nu\underbrace{\int_{\of} \hd^2((\a_j^l \A_j^k-\delta_j^l\delta_j^k) v,_k)\cdot\hd^2 v,_ldxdt}_{I_4}\n\\
	&-\u{\int_0^t\int_{\of} \hd^2(\a_i^j q) \hd^2 v,_j^i dx dt}_{I_5} + \ud\int_{\os} |\hd^2 v|^2(\cdot,t)\ dx+ \frac{\lambda}{2} \int_{\os}|\nabla \hd^2\eta|^2(\cdot,t)\ dx+ g \u{\int_0^t  \int_\of \hd^2\J \hd^2v^3 dx dt}_{I_6}\n\\
	& = \ud\int_{\of} \J|\hd^2 v|^2(\cdot,0)\ dx+\ud\int_{\os} |\hd^2 v|^2(\cdot,0)\ dx  + \frac{\lambda}{2} \int_{\os}|\nabla \hd^2\eta|^2(\cdot,0)\ dx \le C_\lambda E(0)
	\,.
	\label{var1h3}
\end{align}	
Note that we should first use as test function an horizontally convoluted approximation $\hd^2(\rho\ast\hd^2(\rho\ast v))$, and go to the limit as the parameter of approximation converges to zero. The limit process gives the final inequality (\ref{var1h3}). We will not repeat this later on when we do the first time differentiated version of (\ref{var1h3}).

All the terms above can be estimated in a standard way, except $I_6$ and the parts of $I_5$ where $q$ appears without derivative, which present quite formidable obstacles for the global in time existence.

We first explain how $I_1$ can be treated easily. We have
\begin{equation*}
	|I_1|\lesssim \int_0^t  \|\J_t\|_{L^4(\of)} \|\hd^2 v\|_{L^4(\of)} \|\hd^2 v\|_{L^2(\of)}\ ds\n\\
	\lesssim \int_0^t   \|\nabla \tilde v\|_{L^4(\of)}\|\hd^2 v\|_{L^4(\of)} \|\hd^2 v\|_{L^2(\of)}\ ds\,,
		\end{equation*}
		where we used our assumption (\ref{tildeapriori}) for $\J_t$.
		
Sobolev's embeddings then imply
\begin{align*}
	|I_1|\lesssim \int_0^t   \| \tilde v\|_{H^2(\of)}\|\hd^2 v\|_{H^1(\of)} \|\hd^2 v\|_{L^2(\of)}\ ds\,.
\end{align*}
Using (\ref{te}) and (\ref{tildev}), this also implies:
\begin{equation}
	|I_1|\lesssim  \int_0^t \|v\|_{H^2(\of)}\|v\|_{H^3(\of)}^2 ds
	\lesssim  \sup_{[0,t]}(\| v\|_{H^2(\of)})\int_0^t \|v\|_{H^3(\of)}^2 ds
	\lesssim   \sup_{[0,t]} P(N(s)) E(t)\,.\label{I3}
\end{equation}

In a similar way, we have
\begin{equation}
|I_2|+|I_3|+|I_4|\lesssim \sup_{[0,t]} P(N(s)) E(t)\,.\label{I12}
\end{equation}

For $I_6$, we have to make use of the divergence condition in the Stokes problem (\ref{he}) defining $\e$, which shows that $\nabla\J$ is a constant plus a polynomial of degree $\ge 2$ in the gradient of $\e$:
\begin{align}
\J=&\text{det}\nabla(\text{Id}+\e)\n\\
=&1+\div\e+B(\nabla\e)+C(\nabla\e)\n\\
=&1-\frac{1}{|\of|} \int_\g\eta^3 dx_h +B(\nabla\e)+C(\nabla\e)\,,\label{2111.1}
\end{align}
where $B$ and $C$ are bilinear and trilinear forms, with constant coefficients. Moreover, in each product of two or three terms, at least one of the derivatives is an horizontal derivative.

From (\ref{2111.1}) and this structure on horizontal derivatives stated above, we infer
\begin{align}
\|\hd^2 \J\|_{L^2(\of)}= & \|\hd^2 (B(\nabla\e)+C(\nabla\e))\|_{L^2(\of)}\n\\
\lesssim & \|\e\|_{H^3(\of)}\|\hd\e\|_{H^2(\of)} + \|\e\|_{H^3(\of)}^2\|\hd\e\|_{H^2(\of)} \n\\
\lesssim &  N(t)(1+N(t)) \|\hd\e\|_{H^{\frac{3}{2}}(\g)}\label{2111.2}
\end{align}
where we used (\ref{te}) and (\ref{te.2}) in (\ref{2111.2}). From (\ref{2111.2}) we infer
\begin{align}
|I_6|&\lesssim \int_0^t N(t)(1+N(t)) \|\hd\e\|_{H^{\frac{3}{2}}(\g)} \|v\|_{H^2(\of)} dt\n\\ 
&\lesssim \sup_{[0,t]} P(N)\int_0^t \|\hd\e\|_{H^{\frac{3}{2}}(\g)}^2+ \|v\|_{H^2(\of)}^2 dt\n\\
&\lesssim \sup_{[0,t]} P(N)(\int_0^t \|\hd\e\|_{H^{\frac{3}{2}}(\g)}^2 dt+E(t))\,.\label{2111.3}
\end{align}	
Using our fundamental estimate (\ref{inter15}) in (\ref{2111.3}) then yields ($\lambda$ is at least of order $1$)
\begin{align}
	|I_6|
	&\lesssim \sup_{[0,t]} P(N)((1+\sup_{[0,t]} P(N))E(t)+C_\lambda E(0)+E(t))\n\\
		&\lesssim \sup_{[0,t]} P(N)(E(t)+C_\lambda E(0))
	\,.\label{2111.4}
\end{align}

We now turn our attention to the most difficult term to treat, $I_5$. 
\begin{align}
I_5=\u{\int_0^t\int_{\of} \hd^2\a_i^j q \hd^2 v,_j^i dx dt}_{I_{51}} + 2\u{\int_0^t\int_{\of} \hd\a_i^j \hd q \hd^2 v,_j^i dx dt}_{I_{52}}+ \u{\int_0^t\int_{\of} \a_i^j \hd^2 q \hd^2 v,_j^i dx dt}_{I_{53}}\,.
\end{align}
Due to our fundamental estimate (\ref{0611.2}) for $\d\int_0^t \|\hd q\|_{H^1(\of)}^2 dt$, we have in the same way as we proved (\ref{2111.4}) that
\begin{equation}
|I_{52}|\lesssim \sup_{[0,t]} P(N) (E(t)+C_\lambda E(0))\,.\label{I423}
\end{equation}
For $I_{53}$, the situation is slightly more complicated, as the integrand seems like a square (and not a cubic power like the previously estimated term), as $\a$ is close to the identity matrix. However, due to incompressibility, we have
$$\a_i^j \hd^2 v,_j^i=-\hd^2\a_i^j  v,_j^i-2 \hd \a_i^j\hd v,_j^i\,, $$
which shows that the  a priori linear term $\a_i^j \hd^2 v,_j^i$ is in fact a square and this leads to the same type of estimate as earlier:
\begin{equation}
	|I_{53}|\lesssim \sup_{[0,t]} P(N) (E(t)+C_\lambda E(0))\,.\label{I423.bis}
\end{equation}

We now turn our attention to $I_{51}$, that we will need to split in two parts first.

\begin{equation}
I_{51}=\u{\int_0^t (q_e+m(q-q_e))\int_{\of} \hd^2\a_i^j  \hd^2 v,_j^i dx dt}_{I_{511}}+\u{\int_0^t\int_{\of} \hd^2\a_i^j (q-q_e-m(q-q_e)) \hd^2 v,_j^i dx dt}_{I_{512}}
\end{equation}
where the average of $q(\cdot,t)-q_e(\cdot)$ on $\g$ is
$$m(q-q_e)(t)=\frac{1}{|\Gamma|}\int_\g q(\cdot,t)-q_e(\cdot) dx_h\,.$$

Using our $L^\infty$ Poincar\'e-Wirtinger estimate (\ref{qinfty}), we then see that

\begin{align}
|I_{512}|\lesssim & \int_0^t \|\hd^2 \nabla\e\|_{L^2(\of)} ( D (1+P(N))+g N)\|\hd^2\nabla v\|_{L^2(\of)} dt\n\\
\lesssim & \int_0^t \|\hd^2 \nabla\e\|_{L^2(\of)} (D^2 (1+P(N))+g D  N) dt\label{0611.3}\,.
\end{align}

Since $\|\hd^2 \nabla\e\|_{L^2(\of)}\lesssim  N$ by (\ref{te}), we infer from (\ref{0611.3}) that
\begin{align}
	|I_{512}|
	\lesssim & \int_0^t N  D^2 (1+P(N))+g \|\hd^2 \nabla\e\|_{L^2(\of)} D N dt\n\\
	\lesssim& \sup_{[0,t]} P(N) E(t) + g \sup_{[0,t]}N (\int_0^t D^2 dt+\int_0^t \|\hd^2\nabla\e\|_{L^2(\of)}^2 ds)\n\\
	\lesssim& \sup_{[0,t]} P(N) ( E(t) +g\int_0^t \|\hd^2\nabla\e\|_{L^2(\of)}^2 ds)
	  \label{0611.4}\,.
\end{align}

Due to our fundamental estimate (\ref{inter15}), we then finally infer from (\ref{0611.4}) 
\begin{align}
	|I_{512}| 
	\lesssim& \sup_{[0,t]} P(N) ( E(t) +\frac{g}{\lambda} (E(t) (1+\sup_{[0,t]} P(N))+C_\lambda E(0)))\n\\
	\lesssim& \sup_{[0,t]} P(N) (E(t) (1+\sup_{[0,t]} P(N))+C_\lambda E(0))\n\\
	\lesssim& \sup_{[0,t]} P(N) (E(t) (1+C_\lambda E(0))
	  \,.\label{I412}
\end{align}
Note that above we used our assumption that $g$ is small relative to $\lambda$.

To treat $I_{511}$, we notice that due to (\ref{177.2bis}), we have
\begin{align*}
	|I_{511}|
	\lesssim & \int_0^t ( P_\lambda (N)   +  g )\|\hd^2 \nabla\e\|_{L^2(\of)} D dt\n\\
			\lesssim &  (g +\sup_{[0,t]} P_\lambda(N)) \int_0^t {\sqrt{\lambda}}\|\hd^2 \nabla\e\|_{L^2(\of)}^2+\frac{1}{\sqrt{\lambda}} D^2 dt\n\\
			\lesssim & (\frac{g}{\sqrt{\lambda}}+\sup_{[0,t]} P_\lambda(N)) E(t)  +  (g  + \sup_{[0,t]} P_\lambda(N))\int_0^t {\sqrt{\lambda}}\|\hd^2 \nabla\e\|_{L^2(\of)}^2 dt\,.
\end{align*}

Using our fundamental estimate (\ref{inter15}) in this inequality we get:
\begin{align}
	|I_{511}|
	\lesssim & (\frac{g}{\sqrt{\lambda}}+\sup_{[0,t]} P_\lambda(N)) E(t)  +  (g  + \sup_{[0,t]} P_\lambda(N))\frac{\sqrt{\lambda}}{\lambda} (E(t) (1+\sup_{[0,t]} P(N))+C_\lambda E(0))\n\\
	\lesssim & (\frac{g}{\sqrt{\lambda}}+\sup_{[0,t]} P_\lambda (N)) (E(t) +C_\lambda E(0))
	\label{I411}\,.
\end{align}

Using (\ref{I411}), (\ref{I412}), (\ref{I423}), (\ref{I423.bis}), (\ref{2111.4}), (\ref{I12}), (\ref{I3}) in (\ref{var1h3}), we finally obtain:
\begin{align}
	\ud\int_{\Omega} |\hd^2 v|^2(\cdot,t)\ dx + \frac{\lambda}{4} \int_{\os}|\nabla \hd^2\eta|^2(\cdot,t)\ dx + &\nu\int_0^t \int_{\of} |\nabla\hd^2 v|^2 \ dx dt\n\\
	&\lesssim
	(\frac{g}{\sqrt{\lambda}}+\sup_{[0,t]} P_\lambda(N)) (E(t)+C_\lambda E(0)) + C_\lambda E(0)\,.
	\label{var1h3.bis}
\end{align}	

\begin{remark}
	Due to periodicity in the directions $x_1$ and $x_2$, we have $\displaystyle \int_\Omega |v,_{12}|^2 dx=\int_\Omega v_{11} v,_{22}$ by integration by parts, and so summing only over the derivatives of the type $\hd^2$ is not preventing us from estimating the cross second derivative $\d\int_\Omega |v,_{12}|^2 dx$.
\end{remark}

\section{Estimate for the first time-differentiated problem}\label{middle}

We use $-\hd^2 v_t\in H^1_0(\Omega)$ as test function (where $\hd^2$ denote any second order partial derivative $\frac{\partial^2}{\partial x_\alpha\partial x_\alpha}$, for $\alpha=1,2$) in (\ref{var2}) and integrate by parts with respect to the horizontal variable to obtain:
\begin{align*}
	0=&\int_{\of} \hd(\J v_{t})_t\cdot \hd v_t\ dx+ \int_{\of} \hd ((v-\tilde v)_l \a_l^j v,_j)_t \cdot \hd v_t\ dx +\nu\int_{\of} \hd(\a_j^l \A_j^k v,_k)_t\cdot\hd v_t,_l\ dx\n\\
	&
	-\int_{\of} \hd(\a_i^j q)_t \hd {v_t},_j^i dx + \int_{\os}\hd v_{tt}\cdot\hd v_t\ dx+\lambda\int_{\os}\nabla \hd v\cdot\nabla\hd v_t\ dx+g\int_\of\hd\J_t \hd v_t^3 dx\n\\
	=&\int_{\of} \J \hd v_{tt}\cdot \hd v_t\ dx+ \int_{\of} \J_t \hd v_{t}\cdot \hd v_t\  +\hd\J_t  v_{t}\cdot \hd v_t+ \hd\J  v_{tt}\cdot \hd v_t  dx+
	\int_{\of} \hd ((v-\tilde v)_l \a_l^j v,_j)_t \cdot \hd v_t\ dx\n\\
	& +\nu\int_{\of} \hd(\a_j^l \A_j^k v,_k)_t\cdot\hd v_t,_l\ dx
	-\int_{\of} \hd(\a_i^j q)_t \hd {v_t},_j^i dx + \int_{\os}\hd v_{tt}\cdot\hd v_t\ dx+\lambda\int_{\os}\nabla \hd v\cdot\nabla\hd v_t\ dx\n\\
	&+g\int_\of\hd\J_t \hd v_t^3 dx
		\,.
\end{align*}

Integrating in time from $0$ to $t$, and summing over $\alpha=1,2$ (while omitting to write the summation symbols for conciseness) yields
\begin{align}
	&\ud\int_{\of} \J |\hd v_t|^2(\cdot,t)\ dx  \u{-\ud\int_0^t \int_\of \J_t |\hd v_t|^2 dxdt}_{J_1}+\u{\int_0^t \int_{\of} \J_t \hd v_{t}\cdot \hd v_t\  +\hd\J_t  v_{t}\cdot \hd v_t+ \hd\J  v_{tt}\cdot \hd v_t  dx dt}_{J_2}\n\\
	& + g\u{\int_0^t \int_{\of} \hd \J_t \hd v_t^3\ dx dt}_{J_6}+ \u{\int_0^t \int_{\of} \hd ((v-\tilde v)_l \a_l^j v,_j)_t \cdot \hd v_t\ dx dt}_{J_3}+ \nu\int_0^t \int_{\of} |\nabla\hd v_t|^2 \ dx dt\n\\
	&+ \nu\underbrace{\int_{\of} \hd((\a_j^l \A_j^k-\delta_j^l\delta_j^k) {v},_k)_t\cdot\hd {v_t},_ldxdt}_{J_4}
	-\u{\int_0^t\int_{\of} \hd(\a_i^j q)_t \hd {v_t},_j^i dx dt}_{J_5}+\ud\int_{\os} |\hd v_t|^2(\cdot,t)\ dx\n\\
	& + \frac{\lambda}{2} \int_{\os}|\nabla \hd v|^2(\cdot,t)\ dx= \ud\int_{\of} \J|\hd v_t|^2(\cdot,0)\ dx +\ud\int_{\os} |\hd v_t|^2(\cdot,0)\ dx + \frac{\lambda}{2} \int_{\os}|\nabla \hd v|^2(\cdot,0)\ dx \le C_\lambda E(0)
	\,.
	\label{var2h2}
\end{align}

The perturbation terms $J_1$ to $J_4$, as well as $J_6$, can be estimated in the same manner as their counterparts $I_1$ to $I_4$, as well as $I_6$, in Section \ref{highspace}. The reason is that in our framework, for the quantities involved, one time derivative corresponds to one space derivative.
 For the term $J_5$, we have the estimate (\ref{nablaqt}) for $\nabla q_t$ and so any part of $J_5$ where $\nabla q$ or $\nabla q_t$ appears is treated similarly as its counterpart in $I_5$. The difficult terms where either $q$ or $q_t$ appear without any space derivative can be treated similarly as $I_5$ from Section \ref{highspace} due to (\ref{tildev}) and (\ref{inter15}).

This then provides
\begin{align}
	\ud\int_{\Omega} |\hd v_t|^2(\cdot,t)\ dx + \frac{\lambda}{2} \int_{\os}|\nabla \hd v|^2(\cdot,t)\ dx +& \nu\int_0^t \int_{\of} |\nabla\hd v_t|^2 \ dx dt\n\\
	&\lesssim
	(\frac{g}{\sqrt{\lambda}}+\sup_{[0,t]} P_\lambda(N)) (E(t)+C_\lambda E(0)) + C_\lambda E(0)\,.
	\label{var2h2.bis}
\end{align}	

\section{Estimate for the second time-differentiated problem}\label{end}

We formally (see remark below and before (\ref{var3})) use $v_{tt}\in H^1_0(\Omega)$ (which would seem problematic as $v_{tt}$ is a priori only in $L^2(\os)$) as test function in (\ref{var3}) to obtain:
\begin{align*}
	0=&\int_{\of} (\J v_{t})_{tt}\cdot  v_{tt}\ dx+ \int_{\of}  ((v-\tilde v)_l \a_l^j v,_j)_{tt} \cdot v_{tt}\ dx +\nu\int_{\of} (\a_j^l \A_j^k v,_k)_{tt}\cdot v_{tt},_l\ dx\n\\
	&
	-\int_{\of} (\a_i^j q)_{tt}  {v_{tt}},_j^i dx + \int_{\os} v_{ttt}\cdot v_{tt}\ dx+\lambda\int_{\os}\nabla  v_t\cdot\nabla v_{tt}\ dx+g\int_\of \J_{tt} v_{{tt}}^3 dx \n\\
	=&\int_{\of} \J v_{ttt}\cdot  v_{tt}\ dx+ \int_{\of} 2 \J_t  v_{tt}\cdot  v_{tt}\  +\J_{tt}  v_{t}\cdot v_{tt}  dx+
	\int_{\of} ((v-\tilde v)_l \a_l^j v,_j)_{tt} \cdot v_{tt}\ dx +g\int_\of \J_{tt} v_{{tt}}^3 dx\n\\
	& +\nu\int_{\of} (\a_j^l \A_j^k v,_k)_{tt}\cdot v_{tt},_l\ dx
	-\int_{\of} (\a_i^j q)_{tt}  {v_{tt}},_j^i dx + \int_{\os}v_{ttt}\cdot v_{tt}\ dx+\lambda\int_{\os}\nabla  v_t\cdot\nabla v_{tt}\ dx
	\,.
\end{align*}

Integrating in time from $0$ to $t$ yields
\begin{align}
	&\ud\int_{\of} \J | v_{tt}|^2(\cdot,t)\ dx  \u{-\ud\int_0^t \int_\of \J_t | v_{tt}|^2 dxdt}_{K_1}+\u{\int_0^t \int_{\of} 2 \J_t  v_{tt}\cdot  v_{tt}\  +\J_{tt}  v_{t}\cdot v_{tt} dx dt}_{K_2} \n\\
	& + g\u{\int_0^t \int_{\of} \J_{tt}  v_{tt}^3\ dx dt}_{K_6}+ \u{\int_0^t \int_{\of}  ((v-\tilde v)_l \a_l^j v,_j)_{tt} \cdot  v_{tt}\ dx dt}_{K_3}+ \nu\int_0^t \int_{\of} |\nabla v_{tt}|^2 \ dx dt\n\\
	&+ \nu\underbrace{\int_{\of} ((\a_j^l \A_j^k-\delta_j^l\delta_j^k) {v},_k)_{tt}\cdot {v_{tt}},_ldxdt}_{K_4}
	-\u{\int_0^t\int_{\of} (\a_i^j q)_{tt}  {v_{tt}},_j^i dx dt}_{K_5}+\ud\int_{\os} | v_{tt}|^2(\cdot,t)\ dx\n\\
	& + \frac{\lambda}{2} \int_{\os}|\nabla  v_t|^2(\cdot,t)\ dx= \ud\int_{\of} \J| v_{tt}|^2(\cdot,0)\ dx +\ud\int_{\os} | v_{tt}|^2(\cdot,0)\ dx + \frac{\lambda}{2} \int_{\os}|\nabla  v_t|^2(\cdot,0)\ dx \le C_\lambda E(0)
	\,.
	\label{var2h2}
\end{align}	

\begin{remark}
	Note that to justify the inequality (\ref{var2h2}), we can proceed as in \cite{CS2}. We first add a term in $-\kappa\Delta v$ to the linear wave equation (making the problem parabolic-parabolic, and allowing for $\sqrt{\kappa} v_{tt}\in L^2(0,T;H^1_0(\Omega))$) for which the inequality above (with some forcing perturbation terms with $\kappa$, in order to make the compatibility conditions satisfied at time zero) is satisfied, and then pass to the limit as $\kappa\rightarrow 0$, recovering the inequality (\ref{var2h2}) at the limit.
	To be more precise, the inequality which is valid at the limit is when $K_5$ is replaced by its expression obtained below in (\ref{k8}) and (\ref{i42}).
\end{remark}

The perturbation terms $K_1$ to $K_4$, as well as $K_6$, can be estimated in the same manner as their counterparts $I_1$ to $I_4$, as well as $I_6$, in Section \ref{highspace}. The reason is that in our framework, for the quantities involved, one time derivative corresponds to one space derivative. 

The term $K_5$ is where the difficulty appears. The terms where either $q$ or $q_t$ appear without any space derivative can be treated similarly as $I_5$ from Section \ref{highspace} due to (\ref{tildev}) and (\ref{inter15}). The only new difficult term is the part of $K_5$ where $q_{tt}$ appears, due to the lack of good estimate on $q_{tt}$.

\begin{equation}
	\label{k8}
	K_5=\u{\int_0^t\int_{\of} (\a_i^j)_{tt} q v_{tt}^i,_j dxdt}_{K_{51}} + 2\u{\int_0^t\int_{\of} (\a_i^j)_{t} q_t v_{tt}^i,_j  dxdt}_{K_{52}}+ \u{\int_0^t\int_{\of} \a_i^j q_{tt} v_{tt}^i,_j dxdt}_{K_{53}}\,.
\end{equation}

For the same reasons as for $I_5$, we have
\begin{equation}
	|K_{51}|+|K_{52}|\lesssim (\frac{g}{\sqrt{\lambda}}+\sup_{[0,t]} P_\lambda (N))( E(t) + C_\lambda E(0))\,.\label{k7812}
\end{equation}

We now treat $K_{53}$, for which the main challenge is the lack of good estimate for $q_{tt}$. Integrating by parts in time,
\begin{equation*}
	K_{53}=-\int_0^t \int_{\of} \a_i^j {q_{t}} v_{ttt}^i,_j dx dt - \int_0^t \int_{\of} (\a_i^j)_t q_{t} v_{tt}^i,_j dx dt + \bigl[\int_{\of} \a_i^j q_{t} v_{tt}^i,_j dx\bigr]_0^t\,.
\end{equation*}
Using the divergence free condition (\ref{divale}) in the first and last term on the right hand side yields:
\begin{align}
	K_{53}=&3 \u{\int_0^t \int_{\of} (\a_i^j)_t {q_{t}} v_{tt}^i,_j dx dt}_{K_{531}} + 3\u{\int_0^t \int_{\of} (\a_i^j)_{tt} {q_{t}} v_{t}^i,_j dx dt}_{K_{532}}+ \u{\int_0^t \int_{\of} (\a_i^j)_{ttt} {q_{t}} v^i,_j dx dt}_{K_{533}}\n\\
	&- \u{\int_0^t \int_{\of} (\a_i^j)_t q_{t} v_{tt}^i,_j dx dt}_{K_{534}} - 2 \bigl[\int_{\of} (\a_i^j)_t {q_{t}} v_{t}^i,_j dx\bigr]_0^t-\bigl[\int_{\of} (\a_i^j)_{tt} {q_{t}} v^i,_j dx\bigr]_0^t\,.\label{i42}
\end{align}
Using (\ref{qtl4}) and (\ref{tildev}), we obtain in the same way as (\ref{I3}) that
\begin{equation}
	\sum_{i=1}^4|K_{53i}|\lesssim \sup_{[0,t]} P_\lambda (N) E(t)\,.\label{k831to4}
\end{equation}

We now explain how to estimate the last term of (\ref{i42}), as it has some extra difficulty compared to the space-time integrals of (\ref{i42}). The term before is treated the exact same way. The main issue is that the obvious inequality 
\begin{align*}
	\left|\int_{\of} (\a_i^j)_{tt} {q_{t}} v^i,_j dx\right|\lesssim & (\|\nabla\tilde v_t\|_{L^2(\of)}(1+P(N))+\|\nabla v\|^2_{L^4(\of)})\|q_t\|_{L^2(\of)}\|\nabla v\|_{L^\infty(\of)}\,,
\end{align*}
seems to ask too much regularity of $\nabla v$ at each time (our framework has just $v\in L^\infty(0,t;H^2(\of)))$. 
Now, as it turns out, the one space higher regularity $v\in L^\infty(0,t;H^3(\of))$ is controlled independently of time as a consequence of our functional framework and of the properties of the fluid-elastic interaction problem, even though it does not seem obvious a priori. The argument is as follows.
Given that $v$ satisfies for each time the Stokes problem with Dirichlet boundary condition on $\Gamma_{top}$ and Neumann boundary condition on $\Gamma$:
\begin{subequations}
	\label{vh3}
	\begin{alignat}{2}
		-\Delta_{\tilde X}v+\nabla_{\tilde X} (q-q_e)&=-(v_t+(v_j-\tilde v_j)\a_l^j v,_j)-(\nabla_{\tilde X}-\nabla){q_e} \ \  \ &&\text{in} \ \ \of \,,\\
		\div_{\tilde X} v&=0\ \ \ \ &&\text{in} \ \ \of \,,\\
			-\nu \a_j^3 \A_j^k {v^f},_k^i +\a_i^3 (q-q_e)& =-q_e(\a_i^3 - \delta_i^3)  - \lambda (\eta^s-\eta_e),_3^i\,,  &&\text{on} \ \ \g \,,\label{vh3.c}\\
		v &= 0\ \  && \text{on} \ \ \Gamma_{top} \,,
	\end{alignat}
\end{subequations}
where we used (\ref{0211.1}) in (\ref{vh3.c}), we then have by elliptic regularity for (\ref{vh3}) that for each time when the solution remains close to equilibrium
\begin{align}
	\|v\|_{H^3(\of)}+\|q-q_e\|_{H^2(\of)}\lesssim& (1+P(\|\tilde\eta\|_{H^3(\of)}))(\lambda\|(\eta^s-\eta_e),_3\|_{H^{\frac{3}{2}} (\g)} + \sum_{i=1}^3 \|\a_i^3-\delta_i^3\|_{H^{\frac{3}{2}} (\g)})\n\\
	&+ (1+P(\|\tilde\eta\|_{H^3(\of)})) ( \|v_t+(v_i-\tilde v_i)v,_i\|_{H^1(\of)}+\|\a-\text{Id}\|_{H^1(\of)})\n\\
	\lesssim& (1+P(\|\e\|_{H^3(\of)}))(\lambda\|\eta-\eta_e\|_{H^3(\os)}+ \|v_t\|_{H^1(\of)}+\|v\|_{H^2(\of)}\|v,_i\|_{H^1(\of)})\n\\
	&+ (1+P(\|\tilde\eta\|_{H^3(\of)})) (\|\e\|_{H^3(\of)}^2+\|\e\|_{H^3(\of)} )\n\\
	\lesssim& P_\lambda(N(t))\,,\label{higher}
\end{align}
where in the second inequality above we used the fact that each $\a_i^3 - \delta_i^3$ is a sum of derivatives and of products of  derivatives of $\e$. (\ref{higher}) then implies with (\ref{qmultiplier}) that
\begin{align}
	\left|\int_{\of} (\a_i^j)_{tt} {q_{t}} v^i,_j dx\right|\lesssim & N(t) P(N(t)) P_\lambda(N(t))\n\\
\lesssim& N^2(t) P_\lambda(N(t))	\,.\label{197.1}
\end{align}

\begin{remark}
	{\it The higher regularity property (\ref{higher}) was not needed to establish a property such as (\ref{197.1}). Proceeding as in Lemma 6 in (\cite{CS1}), the proof of the Lemma provides that $\|\a_{tt}\|_{L^3(\of)}\le\epsilon+E(t)$, which would have been sufficient.}
\end{remark}

\begin{remark}
	{\it By a similar argument as for (\ref{higher}), we also have the improved regularity}
	\begin{equation*}
		\|v_t\|_{H^2(\of)}+\|q_t\|_{H^1(\of)}\lesssim P_\lambda(N(t))\,,
	\end{equation*}
	{\it even though we will not use this property in this paper.}
\end{remark}

Using (\ref{197.1}), (\ref{k831to4}) in (\ref{i42}) then provides
\begin{equation*}
	|K_{53}|\lesssim \sup_{[0,t]} P_\lambda(N) E(t) + P_\lambda(N(0))E(0)\,,
\end{equation*}
which with our assumption that $N(0)$ is small implies:
\begin{equation*}
	|K_{53}|\lesssim \sup_{[0,t]} P_\lambda(N) E(t) + E(0)\,.
\end{equation*}

Together with (\ref{k7812}) and (\ref{k8}), this provides
\begin{align*}
	|K_{5}|
	\lesssim (\frac{g}{\sqrt{\lambda}}+\sup_{[0,t]} P_\lambda(N))( E(t) + C_\lambda E(0))+C_\lambda E(0) \,.
\end{align*}

This finally implies
\begin{align}
	\ud\int_{\Omega} | v_{tt}|^2(\cdot,t)\ dx + \frac{\lambda}{2} \int_{\os}|\nabla v_t|^2(\cdot,t)\ dx + &\nu\int_0^t \int_{\of} |\nabla v_{tt}|^2 \ dx dt\n\\
	&\lesssim (\frac{g}{\sqrt{\lambda}}+\sup_{[0,t]} P_\lambda(N)) (E(t)+C_\lambda E(0)) + C_\lambda E(0)
	\,.
	\label{var3h1.bis}
\end{align}

\section{Existence for all positive time if the initial data is close enough to equilibrium}
\label{conclusion1}
With (\ref{var3h1.bis}), (\ref{var2h2.bis}) and (\ref{var1h3.bis}), we have
\begin{align*}
\frac 1 2\int_{\Omega} (|v_{tt}|^2+|\hd v_t|^2+|\hd^2 v|^2)(\cdot,t) dx&+\frac{\lambda}{2}\int_{\os} (|\nabla v_t|^2+|\nabla\hd v|^2+|\nabla\hd^2 \eta|^2)(\cdot,t) dx\n\\
&+\nu \int_0^t \int_\of (|\nabla v_{tt}|^2 +|\nabla\hd v_{t}|^2+|\nabla\hd^2 v|^2) dxdt  \lesssim R(t)\,,
\end{align*}	
with 
\begin{equation}
R(t)=(\frac{g}{\sqrt{\lambda}}+\sup_{[0,t]} P_\lambda(N)) (E(t)+C_\lambda E(0)) + C_\lambda E(0)\,.\label{020824.3}
\end{equation}

For obvious reasons, we also could have done the variational estimates (\ref{var2h2.bis}) and (\ref{var1h3.bis}) with a lower order of horizontal derivative, and the estimate would be valid with the same upper bound:
\begin{align}
	&\frac 1 2\int_{\Omega} (|v_{tt}|^2+\sum_{k=0}^1|\hd^k v_t|^2+\sum_{k=0}^2|\hd^k v|^2)(\cdot,t) dx+\frac{\lambda}{2}\int_{\os} (|\nabla v_t|^2+\sum_{k=0}^1|\nabla\hd^k v|^2+\sum_{k=1}^2|\nabla\hd^k \eta|^2)(\cdot,t) dx\n\\
	&+\nu \int_0^t \int_\of (|\nabla v_{tt}|^2 +\sum_{k=0}^1|\nabla\hd^k v_{t}|^2+\sum_{k=0}^2|\nabla\hd^k v|^2) dxdt\n\\
	& +g\int_0^t\int_\of \J v^3 dx dt+ [\int_\os \frac{\lambda}{2}|\nabla\eta|^2 +g\eta^3 dx]_0^t \lesssim R(t)\,.	\label{vartotal.bis}
\end{align}	

As the gravitational terms in the last line above are linear, we rewrite them to make them appear as squares. Also, we remind it is $\eta-\eta_e$ which appears in the norm $N$, with $\eta_e$ given by (\ref{0911.1}). 

We notice that the vertical component part of the variation between $0$ and $t$ in the elastic phase is
\begin{align}
 [\int_\os \frac{\lambda}{2}|\nabla\eta^3|^2 +g\eta^3 dx]_0^t=& [\int_\os \frac{\lambda}{2}|\nabla\eta^3|^2 +g(\eta^3-\eta_e^3) dx]_0^t\n\\
 =& [\int_\os \frac{\lambda}{2}|\nabla(\eta^3-\eta_e^3)|^2 +\lambda \nabla(\eta^3-\eta_e^3)\cdot\nabla\eta_e^3+\frac{\lambda}{2}|\nabla\eta_e^3|^2 +g(\eta^3-\eta_e^3) dx]_0^t\n\\
 =& [\int_\os \frac{\lambda}{2}|\nabla(\eta^3-\eta_e^3)|^2 +\lambda \nabla(\eta^3-\eta_e^3)\cdot\nabla\eta_e^3 +g(\eta^3-\eta_e^3) dx]_0^t\n\\
 =& [\int_\os \frac{\lambda}{2}|\nabla(\eta^3-\eta_e^3)|^2 +\u{(-\lambda \Delta\eta_e^3 +g)}_{=0}(\eta^3-\eta_e^3) dx+\lambda\int_\g (\eta^3-\eta_e^3)\eta_e^3,_3 dx_h]_0^t\n\\
 =& [\int_\os \frac{\lambda}{2}|\nabla(\eta^3-\eta_e^3)|^2 +\frac{gh_s}{2} \int_\g \eta^3 dx_h]_0^t\,.\label{0711.1}
\end{align}

Next, we notice that $(x_h,X^3(X_h^{-1}))$ being a parametrization of $X(\g)$, due to incompressibility, the volume of $X(\of,t)$, and thus of $X(\os,t)$, remains unchanged:
\begin{align}
0=[|X(\os,\cdot)|]_0^t =&[\int_{\Gamma} X^3(X_h^{-1}) dx_h]_0^t\n\\
=&[\int_{X_h^{-1}(\Gamma)} X^3 (X,_1^1X,_2^2-X,_1^2X,_2^1) dy_h]_0^t\n\\
=&[\int_{\Gamma} (h_s+\eta^3) ((1+\eta,_1^1)(1+\eta,_2^2)-\eta,_1^2\eta,_2^1) dx_h]_0^t\,,\label{0711.2}
\end{align} 
where we also used the fact the integral over $X_h^{-1}(\Gamma)$ of any periodic function of period $L$ in the directions $e_1$ and $e_2$ is also the integral over $\Gamma$ of the same function.

Due to horizontal periodicity, integration by parts yields $$\int_\g \eta,_1^1\eta,_2^2 dx_h=\int_\g \eta,_1^2\eta,_2^1 dx_h\,,$$ and so (\ref{0711.2}) becomes
\begin{align}
	0
	=[\int_\g \eta^3 ((1+\eta,_1^1)(1+\eta,_2^2)-\eta,_1^2\eta,_2^1) dx_h]_0^t\,.\label{0711.3}
\end{align} 
The use of (\ref{0711.3}) in (\ref{0711.1}) then yields:
\begin{align}
	[\int_\os \frac{\lambda}{2}|\nabla\eta^3|^2 +g\eta^3 dx]_0^t
	=& [\int_\os \frac{\lambda}{2}|\nabla(\eta^3-\eta_e^3)|^2 -\frac{gh_s}{2} \int_\g \eta^3 (\eta,_1^1+\eta,_2^2+\eta,_1^1\eta,_2^2-\eta,_1^2 \eta,_2^1) dx_h]_0^t\n\\
	=& [\int_\os \frac{\lambda}{2}|\nabla(\eta^3-\eta_e^3)|^2 -\frac{gh_s}{2} \int_\g (\eta^3-\eta_e^3) (\eta,_1^1+\eta,_2^2+\eta,_1^1\eta,_2^2-\eta,_1^2 \eta,_2^1) dx_h]_0^t
	\,,\label{0711.4}
\end{align}
where we used $\eta^3_e=0$ on $\g$ in (\ref{0711.4}). 
Using integration by parts in $\os$ and $\eta=0$ on $\Gamma_B$ for the square part of the integral associated with gravity on $\g$ (and leaving the cubic unchanged), (\ref{0711.4}) becomes
\begin{align}
	[\int_\os \frac{\lambda}{2}|\nabla\eta^3|^2 +g\eta^3 dx]_0^t 
	=& [\int_\os \frac{\lambda}{2}|\nabla(\eta^3-\eta_e^3)|^2 dx]_0^t  -\frac{gh_s}{2}  [\int_\os (\eta^3-\eta_e^3),_3 (\eta,_1^1+\eta,_2^2)dx]_0^t\n\\
	& - \frac{gh_s}{2}[\int_\os (\eta^3-\eta_e^3) (\eta,_{31}^1+\eta,_{32}^2) dx ]_0^t-\frac{gh_s}{2}[\int_\g (\eta^3-\eta^3_e)\text{det}\hd\eta_h dx_h]_0^t\n\\
	=& [\int_\os \frac{\lambda}{2}|\nabla(\eta^3-\eta_e^3)|^2 dx]_0^t  -\frac{gh_s}{2}  [\int_\os (\eta^3-\eta_e^3),_3 (\eta,_1^1+\eta,_2^2)dx]_0^t\n\\
	& +\frac{gh_s}{2} [\int_\os (\eta^3-\eta_e^3),_1 \eta,_{3}^1+ (\eta^3-\eta_e^3),_2\eta,_{3}^2) dx ]_0^t\n\\
	&-\frac{gh_s}{2}[\int_\g (\eta^3-\eta^3_e)\text{det}\hd(\eta)_h dx_h]_0^t\n\\
	=& [\int_\os \frac{\lambda}{2}|\nabla(\eta^3-\eta_e^3)|^2 dx]_0^t  -\frac{gh_s}{2}  [\int_\os (\eta^3-\eta_e^3),_3 ((\eta-\eta_e),_1^1+(\eta-\eta_e),_2^2)dx]_0^t\n\\
	& +\frac{gh_s}{2} [\int_\os (\eta^3-\eta_e^3),_1 (\eta-\eta_e),_{3}^1+ (\eta^3-\eta_e^3),_2(\eta-\eta_e),_{3}^2) dx ]_0^t\n\\
	&-\frac{gh_s}{2}[\int_\g (\eta^3-\eta^3_e)\text{det}\hd(\eta-\eta_e)_h dx_h]_0^t
			\,,\label{0711.5}
\end{align}
where we used in the last equality the fact that $(\eta_e)_h=(0,0)$.

Using Young's inequality, we infer from (\ref{0711.5}) that
\begin{equation*}
	[\int_\os \frac{\lambda}{2}|\nabla\eta^3|^2 +g\eta^3 dx]_0^t \ge \frac{\lambda}{2} \int_\os|\nabla(\eta^3-\eta_e^3)|^2 dx - gh_s \int_\os |\nabla (\eta-\eta_e)|^2 dx  - C_\lambda E(0) -C N(t)^3\,,
\end{equation*}
where we used $N(0)^3\le N(0)^2=E(0)$ for $N(0)\le 1$.

 Adding the gradient of the horizontal components of $\eta$ then yield:
\begin{equation}
	[\int_\os \frac{\lambda}{2}|\nabla\eta|^2 +g\eta^3 dx]_0^t \ge \frac{\lambda}{2} \int_\os|\nabla(\eta-\eta_e)|^2 dx - gh_s \int_\os |\nabla (\eta-\eta_e)|^2 dx  - C_\lambda E(0) -C N(t)^3\,.\label{0711.6}
\end{equation}

We now study the space-time integral associated with the gravitational term in $\of$ in (\ref{vartotal.bis}):

\begin{equation}
	\int_0^t\int_\of \J v^3 dx dt=	\int_0^t\int_\of \J u^3(\tilde X) dx dt
	=	\int_0^t\int_{\tilde X(\of)} u^3 dX dt\,.\label{2311.1}
\end{equation}
Now, since the Lagrangian flow map $X$ also produces the same fluid domain as $\tilde X$ (namely the same infinite periodic fluid domain), we have that the integrals of any periodic function of period $L$ in the directions $e_1$ and $e_2$ over $\tilde X(\of)$ and $X(\of)$ are equal. From (\ref{2311.1}) we then infer:
\begin{equation*}
	\int_0^t\int_\of \J v^3 dx dt=	\int_0^t\int_{X(\of)} u^3 dX dt
	=	\int_0^t\int_{\of} \text{det}\nabla X\ u^3(X) dx dt
	=	\int_0^t\int_{\of} \text{det}\nabla X\ X^3_t dx dt\,.
\end{equation*}
Now, since $(\text{det}\nabla X)_t=0$ for the Lagrangian map of an incompressible fluid, we infer from above that
\begin{equation*}
	\int_0^t\int_\of \J v^3 dx dt
	=[\int_\of \text{det}\nabla X\ X^3 dx]_0^t=[\int_{X(\of,t)} x_3 dx]_0^t\,.
\end{equation*}

Next,  due to $|X(\of,t)|$ remaining constant in time by incompressibility, we infer from above that
\begin{equation*}
	\int_0^t\int_\of \J v^3 dx dt
	=[\int_{X(\of,t)} (x_3-h_s) dx]_0^t=\ud[\int_{X(\of,t)} (x_3-h_s)^2,_3 dx]_0^t\,.
\end{equation*}

Integrating by parts in $\Omega^f(t)$, we then infer
\begin{equation*}
	\int_0^t\int_\of \J v^3 dx dt
	=\ud[\int_{X(\Gamma,t)} (x_3-h_s)^2 n_3 dS+ \int_{\Gamma_{top}} (h-h_s)^2 dx_h]_0^t
		=\ud[\int_{X(\Gamma,t)} (x_3-h_s)^2 n_3 dS]_0^t
	\,.
\end{equation*}
Going back to the reference interface, this yields:
\begin{align}
\int_0^t\int_\of \J v^3 dx dt
=&-\ud[\int_{\Gamma} (h_s+\eta^3-h_s)^2 |X,_1\times X,_2| \frac{a_3^3}{|a_i^3|} dx_h]_0^t\n\\
=&-\ud[\int_{\Gamma} (\eta^3-\eta_e^3)^2 {(1+\eta,_1^1+\eta,_2^2+\eta,_1^1\eta,_2^2-\eta,_1^2\eta,_2^1)}  dx_h]_0^t\n\\
=&-\ud[\int_{\Gamma} (\eta^3-\eta_e^3)^2 {(1+(\eta-\eta_e),_1^1+(\eta-\eta_e),_2^2+\det\hd(\eta-\eta_e)^h)}  dx_h]_0^t
\,,\label{0811.1}
\end{align}
where we have used $\eta_e=0$ on $\g$ in (\ref{0811.1}). By the fundamental theorem of calculus in $\os$, and Cauchy-Schwarz used with $\d\eta^3-\eta_e^3=\int_0^{h_s} (\eta^3-\eta_e^3),_3 dx_3$ on $\g$, we infer from (\ref{0811.1}) that
\begin{align}
	g\int_0^t\int_\of \J v^3 dx dt
	\ge &-\frac{1}{2} gh_s\int_\os (\eta^3,_3-\eta_e^3,_3)^2 dx-C E(0)- C (N^3(t)+N^4(t))\,.\label{0811.2}
\end{align}
We then infer from (\ref{0711.6}) and (\ref{0811.2}) that
\begin{equation*}
	g \int_0^t\int_\of \J v^3 dx dt+[\int_\os \frac{\lambda}{2}|\nabla\eta|^2 +g\eta^3 dx]_0^t\ge (\frac{\lambda}{2}-\frac{3}{2}gh_s) \int_\os |\nabla(\eta-\eta_e)|^2 dx - C_\lambda E(0) -C N^2(t) P(N(t))\,.
\end{equation*}
Using our assumption (\ref{lambdalarge}) on $\lambda$ compared to $g$, this implies
\begin{equation}
	g \int_0^t\int_\of \J v^3 dx dt+[\int_\os \frac{\lambda}{2}|\nabla\eta|^2 +g\eta^3 dx]_0^t\ge \frac{\lambda}{4}\int_\os |\nabla(\eta-\eta_e)|^2 dx - C_\lambda E(0) -C E(t) P(N(t))\,.\label{0811.3}
\end{equation}

Using (\ref{0811.3}) in (\ref{vartotal.bis}) then finally provides
\begin{align}
	&\frac 1 2\int_{\Omega} (|v_{tt}|^2+\sum_{k=0}^1|\hd^k v_t|^2+\sum_{k=0}^2|\hd^k v|^2)(\cdot,t) dx+\frac{\lambda}{2}\int_{\os} (|\nabla v_t|^2+\sum_{k=0}^1|\nabla\hd^k v|^2+\sum_{k=1}^2|\nabla\hd^k \eta|^2)(\cdot,t) dx\n\\
	&+\frac{\lambda}{4}\int_\os |\nabla(\eta-\eta_e)|^2 dx+\nu \int_0^t \int_\of (|\nabla v_{tt}|^2 +\sum_{k=0}^1|\nabla\hd^k v_{t}|^2+\sum_{k=0}^2|\nabla\hd^k v|^2) dxdt\lesssim R(t)\,.	\label{vartotal}
\end{align}

We now remind why the left hand side of (\ref{vartotal}) controls $E(t)$, due to similar arguments as in \cite{CS1}, \cite{CS2}.
By elliptic regularity for the elliptic system with Dirichlet boundary conditions,
 \begin{subequations}
	\label{wave}
	\begin{alignat}{2}
		\Delta(\eta-\eta_e) &=\frac{v_t}{\lambda} \ \  \ &&\text{in} \ \ \os \,,\\
		\eta-\eta_e &= \eta-\eta_e|_\g \ \  &&\text{on} \ \ \g \,,\\
		\eta-\eta_e &= 0\ \  && \text{on} \ \ \Gamma_{B} \,,
	\end{alignat}
\end{subequations}
we have
	\begin{align*}
		\|\eta-\eta_e\|^2_{H^3(\os)}\lesssim & \frac{\|v_t\|^2_{H^1(\os)}}{\lambda^2}+\|\hd^2\eta\|^2_{H^{\ud}(\g)}+\|\eta-\eta_e\|^2_{H^{\ud}(\g)}\,,
	\end{align*}
	where $\hd$ again stands for the full horizontal gradient (and we remind $\eta_e=0$ on $\g$).
	
Due to the properties of the trace, and Poincar\'e's\ inequality, we infer successively from this that
\begin{align}
	\|\eta-\eta_e\|^2_{H^3(\os)}\lesssim & \frac{\|v_t\|^2_{H^1(\os)}}{\lambda^2}+\|\hd^2\eta\|^2_{H^1(\os)}+\|\eta-\eta_e\|^2_{H^1(\os)}\n\\
	\lesssim & \frac{\|\nabla v_t\|^2_{L^2(\os)}}{\lambda^2}+\|\nabla\hd^2\eta\|^2_{L^2(\os)}+\|\nabla(\eta-\eta_e)\|^2_{L^2(\os)}\n\\
	\lesssim & R(t)(1+\frac{1}{\lambda^3})\n\\
	\lesssim & R(t)\,,\label{207.1}
\end{align}
due to (\ref{vartotal}) for the third inequality and $\lambda\ge 1$
for the fourth. By the exact same type of arguments, we also have
\begin{equation}
\|v\|^2_{H^2(\os)}\lesssim R(t)\,.\label{207.2}
\end{equation}

Similarly, elliptic regularity on
\begin{subequations}
	\label{nsd}
	\begin{alignat}{2}
		-\Delta_{\tilde X}v+\nabla_{\tilde X} (q+gx_3)&=-(v_t+(v_i-\tilde v_i)\A_i^lv,_l)+g (-\nabla+\nabla_{\tilde X}) x_3  \ \  \ &&\text{in} \ \ \of \,,\\
		\div_{\tilde X} v&=0\ \ \ \ &&\text{in} \ \ \of \,,\\
		v &= v^s  \ \  &&\text{on} \ \ \g \,,\\
		v &= 0\ \  && \text{on} \ \ \Gamma_{top} \,,
	\end{alignat}
\end{subequations}
yields
\begin{align*}
\|v\|^2_{H^2(\of)}+\|\nabla(q+gx_3)\|^2_{L^2(\of)}\lesssim& \|v_t+(v_i-\tilde v_i)\A_i^l v,_l\|^2_{L^2(\of)}+ \|v\|^2_{H^{\frac 3 2} (\g)}+g^2\|\A_i^3-\delta_i^3\|_{L^2(\of)}^2\n\\
\lesssim& \|v_t\|^2_{L^2(\of)}+\|v-\tilde v\|^2_{L^4(\of)} \|\nabla v\|^2_{L^4(\of)} (1+P(N(t)))+   \|v\|^2_{H^2 (\os)} \n\\
&+g^2 \|\hd\e\|_{L^2(\os)}^2(1+P(N(t)))\n\\
\lesssim& \|v_t\|^2_{L^2(\of)}+\|v-\tilde v\|^2_{H^1(\of)} \|\nabla v\|^2_{H^1(\of)}(1+P(N(t)))+  \| v\|^2_{H^2 (\os)}\n\\
&+g^2 \|\hd\e\|_{H^2(\of)}^2 (1+P(N))\,,
\end{align*}
where we used the trace theorem from $\os$ in the second inequality above. This now becomes:
\begin{align}
	\|v\|^2_{H^2(\of)}+\|\nabla(q+gx_3)\|^2_{L^2(\of)}
\lesssim& \|v_t\|^2_{L^2(\of)}+N(t)^4 (1+P(N(t)))+\| v\|^2_{H^2 (\os)}\n\\
&+g^2 \|\hd(\eta-\eta_e)\|_{H^2(\os)}^2 (1+P(N))\,,\label{020824.1}
\end{align}
where we used (\ref{te.2}) in the last line of (\ref{020824.1}). Therefore, (\ref{207.2}), (\ref{vartotal}) and (\ref{207.1}) allow us to infer from (\ref{020824.1}) that
\begin{align}
	\|v\|^2_{H^2(\of)}+\|\nabla(q+gx_3)\|^2_{L^2(\of)}\lesssim& R(t) +E(t) P(N(t)) +  R(t)\n\\
	&+g^2 R(t) (1+P(N))\n\\
	\lesssim& R(t)\label{207.3}\,,
\end{align}
where we used the stability of the generic form of $R(t)$ by addition and multiplication by $(1+P(N(t)))$ (as we can see from (\ref{020824.3})) in the final inequality.
By the same type of arguments, and our fundamental property (\ref{inter15}) about the $L^2$ in time integrability of $\hd\e$, we obtain similarly:
\begin{equation}
	 \int_0^t \|v\|^2_{H^3(\of)} dt \lesssim R(t)\,.\label{1407.4}
\end{equation}

By working on the variational formulation associated to the time differentiated version of (\ref{nsd}), we can prove similarly that
\begin{equation}
	\|v_t\|^2_{H^1(\of)} +\int_0^t \|v_t\|^2_{H^2(\of)} dt \lesssim R(t)\,.\label{207.4}
\end{equation}
With (\ref{207.1}), (\ref{207.2}), (\ref{207.3}), (\ref{1407.4}) and (\ref{207.4}), we infer from (\ref{vartotal}) that
\begin{equation}
	E(t)\le C_1 R(t) = C_1 (\frac{g}{\sqrt{\lambda}}+\sup_{[0,t]} P_\lambda(N)) ( E(t)+C_\lambda E(0)) +C_1 C_\lambda E(0) \,,\label{18.3}
\end{equation}
where $C_1>0$ finite (and independent of $\lambda$ and $t$) is now fixed.
Due to Young's inequality we can assume our polynomial $P(N)$ under the form:
$$C_1 P_\lambda (N)=C_{2,\lambda} N+C_{2,\lambda} N^{2n}\,,$$ for some $n\ge 2$ integer. Again due to Young's inequality, this implies
\begin{equation}
C_1 P_\lambda (N)\le \frac{\epsilon}{2}+ C_{\epsilon,\lambda} N^{2n}\,.\label{18.4}
\end{equation}
Since $\displaystyle E(t)\ge \sup_{[0,t]} N^2$, we infer from (\ref{18.4}) that
\begin{equation}
	C_1 \sup_{[0,t]} P_\lambda (N)\le \frac{\epsilon}{2}+C_{\epsilon,\lambda} E^{n}\,.\label{18.5}
\end{equation}

Therefore, we infer from (\ref{18.3}) and(\ref{18.5}) that for $\lambda$ large enough so that $$C_1\frac{g}{\sqrt{\lambda}}\le \frac{\epsilon}{2}\,,$$ we have
\begin{align}
E(t)&\le (\epsilon+C_{\epsilon,\lambda} E^n) (E(t)+C_\lambda E(0))+ C_1C_\lambda E(0)\n\\
&\le \epsilon E(t)+ D_{\lambda,\epsilon} E(0)+  D_{\lambda,\epsilon} E^{n+1}(t) + D_{\lambda,\epsilon} E(0) E(t)^n
\,,\label{28.1}
\end{align}
where $D_{\lambda,\epsilon}>0$ depends on $\epsilon$ and $\lambda$.

Again, using Young's inequality on the fourth term of the right hand side of (\ref{28.1}), we infer
\begin{align}
	E(t)&\le \epsilon E(t)+ \tilde D_{\lambda,\epsilon} E(0)+ \tilde D_{\lambda,\epsilon} E^{n+1}(t) + \tilde D_{\lambda,\epsilon} E(0)^n\n\\
	&\le \epsilon E(t)+ 2\tilde D_{\lambda,\epsilon} E(0)+ \tilde D_{\lambda,\epsilon} E^{n+1}(t)\,,\label{28.4}
\end{align}
if we assume our initial data small enough. We can also assume $D_{\lambda,\epsilon}\ge 1$ (if not, any larger value for this number still makes (\ref{28.4}) valid).

We first choose $0<\epsilon<\ud$, so that (\ref{28.4}) implies
 that for all time of existence, so long as the solution exists and $N(t)\le \frac{\lambda}{2C_{0,\lambda}}$ on $[0,T]$ (where $C_{0,\lambda}$ is defined in (\ref{smalldata})):
\begin{align}
	E(t)\le 4\tilde D_{\lambda,\epsilon} E(0)+ 2\tilde D_{\lambda,\epsilon} E^{n+1}(t)\,.\label{28.6}
\end{align}
The variation of the polynomial $$f(x)= 4\tilde D_{\lambda,\epsilon} E(0)+ 2\tilde D_{\lambda,\epsilon} x^{n+1}-x\,,$$
show that $f$ decreases from $f(0)= 4\tilde D_{\lambda,\epsilon} E(0)>0$ (we can assume $E(0)=N(0)^2>0$, otherwise we have the initial data corresponding to the canonical equilibrium) to $f(x_0)$ with $2 \tilde D_{\lambda,\epsilon} (n+1) x_0^{n}=1$, and then increases from $f(x_0)$ to $\infty$.

If we impose our initial data to be small enough in order to satisfy
\begin{equation*}
 4 \tilde D_{\lambda,\epsilon} E(0)+ 2 \tilde D_{\lambda,\epsilon} (8 \tilde D_{\lambda,\epsilon} E(0))^{n+1}-8 \tilde D_{\lambda,\epsilon} E(0)= -4 \tilde D_{\lambda,\epsilon} E(0)+ 2 \tilde D_{\lambda,\epsilon} (8\tilde D_{\lambda,\epsilon})^{n+1}E^{n+1}(0)< 0\,,
\end{equation*}
this means by definition that we impose
\begin{equation}
f(8 \tilde D_{\lambda,\epsilon} E(0))< 0\,.\label{28.7}
\end{equation}

We moreover impose our initial data to satisfy
\begin{equation}
8 \tilde D_{\lambda,\epsilon} E(0)<x_0\,,\label{28.8}
\end{equation}
which implies $8 \tilde D_{\lambda,\epsilon} E(0)\in (z_0,z_1)$, where $z_0$ and $z_1$ are the two zeros of $f$ on $[0,\infty)$.

Therefore, $f\ge 0$ on $[0,z_0]$ and $[z_1,\infty)$, with
\begin{equation}
	0<z_0<8 \tilde D_{\lambda,\epsilon} E(0)<z_1\,.\label{order1}
\end{equation}

We have already established in (\ref{28.6}) that so long as the solution exists and the small data assumption (\ref{smalldata}) is satisfied on $[0,T(E(0))]$:
\begin{equation}
f(E(t))\ge 0\,.\label{28.9}
\end{equation}

Due to (\ref{order1}), 
\begin{enumerate}
\item either for all $t\in [0,T(E(0))]$, $E(t)\ge z_1>8 \tilde D_{\lambda,\epsilon} E(0)$,

\item or for all $t\in [0,T(E(0))]$, $E(t)\le z_0<8 \tilde D_{\lambda,\epsilon} E(0)$.
\end{enumerate}

Due to $E(0)<8 E(0)\le 8  \tilde D_{\lambda,\epsilon} E(0)$, we see that it is proposition 2 above which is true.

This implies from $f(0)>0$ and (\ref{28.7}) that so long as the solution exists and $N(t)\le \frac{\lambda}{2C_{0,\lambda}}$ on $[0,T]$:
\begin{equation}
	E(t)< 8 \tilde D_{\lambda,\epsilon} E(0)=8 \tilde D_{\lambda,\epsilon} N(0)^2 \,.\label{28.10}
\end{equation}
If we moreover impose that
\begin{equation}
8 \tilde D_{\lambda,\epsilon} E(0) < \frac{\lambda^2}{4{C_{0,\lambda}}^2}\,,\label{28.11}
\end{equation}
we then have as a consequence of (\ref{28.10}) that on $[0,T(E(0))]$
\begin{equation}
N(t)\le \sqrt{E(t)}\le\frac{\lambda}{2 C_{0,\lambda}}\,,
\end{equation}
and so the condition (\ref{smalldata}) $N(t)\le \frac{\lambda}{2C_{0,\lambda}}$ needed on $[0,T(E(0))]$ is automatically satisfied. We also remind that we picked $C_{0,\lambda}$ large enough so that (\ref{smalldata}) being satisfied implies that (\ref{tildeapriori}) is also satisfied. Therefore, we infer that for an initial data satisfying all the conditions of smallness of this Section, as the solution exists on $[0,T(E(0))]$, (\ref{28.10}) is satisfied on $[0,T(E(0))]$. As it can then be extended (local in time existence) from $[T(E(0)),T(E(0))+T(8\tilde D_{\lambda,\epsilon} E(0))]$, we then have that as the solution exists on $[0,T(E(0))+T({8 \tilde D_{\lambda,\epsilon}} E(0))]$, it satisfies (\ref{28.10}) on $[0,T(E(0))+T( {8 \tilde D_{\lambda,\epsilon}}E(0))]$. By induction, it satisfies (\ref{28.10}) on all intervals $[0,T(E(0))+nT({8 \tilde D_{\lambda,\epsilon}} E(0))]$, and so on $[0,\infty)$.

This finishes the proof of Theorem 1.\qed

\section{Convergence in large time towards a flat interface solution as defined in \ref{flat}}
\label{conclusion_asymptotic}

In this Section, we assume our data close enough to equilibrium initially to satisfy Theorem \ref{theorem_main}. In this Section, we will not need to keep track of the dependence on $\lambda$ in estimates.

\subsection{Convergence as $t\rightarrow\infty$ towards equilibrium in the fluid and at the interface}

Our starting point will be that since the data satisfies Theorem \ref{theorem_main}, we have
\begin{equation}
\int_0^\infty \|v\|^2_{H^2(\of)}+ \|v_t\|^2_{H^2(\of)} dt<\infty\,,\label{cv0}
\end{equation}
which by Cauchy-Schwarz provides $\ddt\|v\|^2_{H^2(\of)}\in L^1(0,\infty)$. This in turn implies the existence of $l$ finite such that
\begin{equation*}
l=\lim_{t\rightarrow\infty}\|v\|^2_{H^2(\of)}\,,
\end{equation*}
which together with $v\in L^2(0,\infty;H^2(\of))$ imply $l=0$, and so the fluid velocity converges to zero:
\begin{equation}
	0=\lim_{t\rightarrow\infty}\|v\|^2_{H^2(\of)}\,.\label{cv1}
\end{equation}
The same arguments allow the same conclusion for the acceleration in the fluid phase:
\begin{equation}
	0=\lim_{t\rightarrow\infty}\|v_t\|^2_{H^1(\of)}\,.\label{cv1.bis}
\end{equation}

For the solid phase, we first conclude about the behaviour of the moving interface, and some of the interior derivative.

From Theorem \ref{theorem_main}, and our fundamental relation (\ref{inter15}), we have
\begin{equation}
\int_0^\infty \|\hd\tilde\eta\|^2_{H^2(\of)}+ \|\tilde v\|^2_{H^3(\of)} dt<\infty\,,\label{cv1.ter}
\end{equation}

we have for the same reason as (\ref{cv1}):
\begin{equation*}
	0=\lim_{t\rightarrow\infty}\|\hd\tilde\eta\|^2_{H^2(\of)}\,.
\end{equation*}
This then implies by the properties of the trace that
\begin{equation}
	0=\lim_{t\rightarrow\infty}\|\hd\tilde\eta\|^2_{H^{\frac 3 2}(\g)}=\lim_{t\rightarrow\infty}\|\hd\eta\|^2_{H^{\frac 3 2}(\g)} \,.\label{cv2}
\end{equation}
Due to volume conservation, since $|\Omega^s(t)|=h_e L^2$, we have for all $t$ the existence of $x_h(t)$ such that
\begin{equation}
\eta^3(x_h(t),h_s,t)=h_e-h_s\,.\label{cv3}
\end{equation}
With (\ref{cv2}) this in turn implies by a Poincar\'e-Wirtinger type inequality that
\begin{equation*}
	0=\lim_{t\rightarrow\infty}\|\eta^3-(h_e-h_s)\|^2_{L^2(\g)} \,,
\end{equation*}
and thus that
\begin{equation}
	0=\lim_{t\rightarrow\infty}\|\eta^3-(h_e-h_s)\|^2_{H^{\frac{5}{2}}(\g)}\,.\label{cv4}
\end{equation}
Obtaining the same type of result for the horizontal component will require more work. 

From the continuity of stress (\ref{systemhorizontal}) along $\g$,
our global in time existence in Theorem \ref{theorem_main} and our higher regularity estimate (\ref{higher}) for the pressure, we then have for the horizontal component
\begin{align}
\int_0^\infty \|(\eta,_3^s)^h\|_{L^2(\g)}^2\lesssim & \int_0^\infty \sum_{\alpha=1}^2 \|\a_\alpha^3\|_{L^2(\g)}^2 dt +\int_0^\infty \|\nabla v^f\|_{L^2(\g)}^2 dt\n\\
\lesssim & \int_0^\infty \|\hd\tilde\eta\|_{L^2(\g)}^2 dt +\int_0^\infty D(t)^2 dt\n\\
<&\infty\,,\label{cv5}
\end{align}
due to (\ref{cv1.ter}) and our definition of the dissipative energy.

We now consider the horizontal component of the linear wave equation, take the scalar product with  $x_3  \eta,_3^h$ and integrate in space-time.

\begin{align*}
	0=&2\int_0^t \int_{\os} v_{t}^h \cdot\eta,_3^h x_3   dx dt -2\lambda \int_0^t\int_{\os} \Delta \eta^h\cdot \eta,_3^h x_3  dx dt\n\\
	=&-2\int_0^t \int_{\os} v^h\cdot v,_3^h x_3   dx dt +2[\int_{\os} v^h\cdot \eta,_3^h x_3 (\cdot)dx]_0^t\n\\
	& +2\lambda\int_0^t\int_{\os} \nabla \eta^h\cdot\nabla \eta,_3^h x_3  +\eta,_3^h\cdot \eta,_3^h  dx dt-2\lambda\int_0^t\int_{\p\os} \eta,_3^h\cdot\eta,_3^h x_3  N_3 dx_h dt\n\\
	=&-\int_0^t \int_{\p\os}|v^h|^2 x_3  N_3 dx_h dt +\int_0^t \int_{\os} |v^h|^2    dx dt +2[\int_{\os} v^h\cdot \eta,_3^h x_3 (\cdot)dx]_0^t\n\\
	& +\lambda\int_0^t \int_{\p\os}|\nabla \eta^h|^2 x_3  N_3 dx_h dt +\lambda\int_0^t\int_{\os} -|\nabla \eta^h|^2   +2 \eta,_3^h\cdot \eta,_3^h   dx dt\n\\
	&-2\lambda\int_0^t\int_{\p\os} \eta,_3^h\cdot \eta,_3^h x_3   N_3 dx_h dt\,.
\end{align*}

Since $x_3=0$ on $\Gamma_B$, and $x_3=h_s$ on $\Gamma$, the previous identity becomes:

\begin{align}
	0=&-h_s \int_0^t  \int_{\g}|v^h|^2 +\lambda|\eta,_3^h|^2-\lambda|\hd\eta^h|^2  dx_h dt +\int_0^t   \int_{\os} |v^h|^2+\lambda|\eta,_3^h|^2-\lambda|\hd\eta^h|^2dx dt\n\\
	& +2\left[\int_{\os} v^h\cdot \eta,_3^h x_3 dx\right]_0^t\,.\label{cv6}
\end{align}
Due to (\ref{cv1.ter}) and (\ref{cv5}) and the fact $N(t)$ is bounded on $[0,\infty)$, we infer from (\ref{cv6}) that there exists $C_1>0$ finite such that
\begin{equation}
\forall t\ge 0\,,\ \ \left|\int_0^t  \int_{\os} |v^h|^2+\lambda|\eta,_3^h|^2-\lambda|\hd\eta^h|^2dx dt\right|\le C_1\,.\label{cv7}
\end{equation}

We still take the horizontal component of the wave equation, but now take the scalar product with $\eta^h$:
 \begin{align}
 0=&\int_0^t\int_\os v_t^h\cdot\eta^h dx dt+\lambda\int_0^t\int_\os |\nabla\eta^h|^2 dxdt -\lambda\int_0^t\int_\g \eta,_3^h\cdot \eta^h dx_h dt\n\\
 =&-\int_0^t\int_\os |v^h|^2dx dt+\left[\int_\os v^h\cdot\eta^h dx\right]_0^t+\lambda\int_0^t\int_\os |\nabla\eta^h|^2 dxdt -\lambda\int_0^t\int_\g \eta,_3^h\cdot \eta^h dx_hdt\,.\label{cv8}
 \end{align}

Using (\ref{cv8}) in (\ref{cv7}) then provides:
\begin{equation*}
	\forall t\ge 0\,,\ \ \left|\lambda\int_0^t\int_\os 2|\eta,_3^h|^2 dx dt-\lambda\int_0^t\int_\g \eta,_3^h\cdot \eta^h dx_hdt  +\left[\int_\os v^h\cdot\eta^h dx\right]_0^t\right|\le C_1\,,
\end{equation*}
which due to $N$ being bounded on $[0,\infty)$ provides the existence of $C_2>0$ finite such that
\begin{equation}
	\forall t\ge 0\,,\ \ \left|\int_0^t\int_\os 2|\eta,_3^h|^2 dx dt-\int_0^t\int_\g \eta,_3^h\cdot \eta^h dx_hdt \right|\le C_2\,.\label{cv9}
\end{equation}

From $\eta=0$ on $\Gamma_B$, we have by the fundamental theorem of calculus and Cauchy-Schwarz that
\begin{align}
\int_\g |\eta^h|^2 dx_h=&\int_{[0,L]^2} \left|\int_0^{h_s} \eta,_3^h dx_3\right|^2 dx_h\n\\
\le& h_s\int_{[0,L]^2} \int_0^{h_s} |\eta,_3^h|^2 dx_3 dx_h=h_s\int_\os|\eta,_3^h|^2 dx \,.\label{cv10}
\end{align} 

Therefore, (\ref{cv9}) and (\ref{cv10}) imply
\begin{equation*}
	\forall t\ge 0\,,\ \ \frac{2}{h_s} \int_0^t\int_\g |\eta^h|^2 dx_h dt-\int_0^t\int_\g \eta,_3^h\cdot \eta^h dx_hdt \le C_2\,.
\end{equation*}
By Young's inequality, this implies
\begin{equation*}
	\forall t\ge 0\,,\ \ \frac{2}{h_s} \int_0^t\int_\g |\eta^h|^2 dx_h dt\le \int_0^t\int_\g \frac{h_s}{2}|\eta,_3^h|^2 + \frac{1}{2h_s}|\eta^h|^2 dx_hdt + C_2\,,
\end{equation*}
which with (\ref{cv5}) finally yields
\begin{equation}
\int_0^\infty \int_\g |\eta^h|^2 dx_h dt<\infty\,.\label{cv11}
\end{equation}
By definition of the dissipative energy, we also have that
\begin{equation}
	\int_0^\infty \int_\g |v^h|^2 dx_h dt<\infty\,.\label{cv12}
\end{equation}
In a way similar as we obtained (\ref{cv1}), the inequalities (\ref{cv11}) and (\ref{cv12}) yield:
\begin{equation*}
	0=\lim_{t\rightarrow\infty}\|\eta^h\|^2_{L^2(\g)} \,,
\end{equation*}
and thus, with (\ref{cv2}) we infer that
\begin{equation}
	0=\lim_{t\rightarrow\infty}\|\eta^h\|^2_{H^{\frac{5}{2}}(\g)}\,.\label{cv13}
\end{equation}

The limits (\ref{cv4}) and (\ref{cv13}) show that the interface $\Gamma(t)$ converges towards the flat interface $\Gamma_e$ in $H^{\frac 5 2}(\g)$ as $t\rightarrow\infty$.
 \subsection{Convergence in $\os$}

 \subsubsection{Convergence of the horizontal displacement and its first and second time derivatives to zero}

Using (\ref{cv5}) and (\ref{cv11}) in (\ref{cv9}), we then infer
\begin{equation}
	\int_0^\infty 	\|\eta,_3^h\|^2_{L^2(\os)} dt<\infty\,.\label{cv15}
\end{equation}

When taking one time derivative of the continuity of stress (\ref{systemhorizontal}) on $\g$, we similarly have
\begin{align*}
	\|v^s,_3^h\|_{L^2(\g)}\lesssim& \|\nu (\a_j^3\A_j^k v,_k^h)_t - q_t a_h^3-q (a_h^3)_t\|_{L^2(\g)}\n\\
	\lesssim & \|\nabla v_t^f\|_{L^2(\g)} + \|\nabla \tilde v\|_{L^4(\g)}\|\nabla v^f\|_{L^4(\g)}  + \|q_t\|_{L^2 (\g)}\|\hd\tilde \eta\|_{L^\infty(\g)}+\|q\|_{L^2 (\g)}\|\hd v\|_{L^\infty(\g)} \n\\
	\lesssim & \| v_t\|_{H^2(\of)} + \|v\|_{H^2(\of)}^2+\|\hd\e\|_{H^2(\of)}+ D(t) \|\hd\e\|_{H^2(\of)}+ \|\hd v\|_{H^2(\of)}
	\n\\
	\lesssim & \| v_t\|_{H^2(\of)} + \|v\|_{H^2(\of)}N(t)+\|\hd\e\|_{H^2(\of)}+ D(t) N(t)+ \|\hd v\|_{H^2(\of)}
	\,,
\end{align*}
where we used our pressure control (\ref{177.2bis}) and (\ref{177.2ter}) for $q$ and $q_t$, and (\ref{te}). We then infer
\begin{equation*}
	\|v^s,_3^h\|_{L^2(\g)}\lesssim D(t)+\|\hd\e\|_{H^2(\of)}\,.
\end{equation*}

This then implies using (\ref{cv1.ter}) and the definition of the dissipative energy that
\begin{equation}
	\int_0^\infty 	\|v^s,_3^h\|^2_{L^2(\g)} dt<\infty\,.\label{cv16}
\end{equation}

We then infer from (\ref{cv16}), in a manner similar as (\ref{cv9}), that we have the existence of $C_3>0$ finite such that
\begin{equation}
	\forall t\ge 0\,,\ \ \left|\int_0^t\int_\os 2|v,_3^h|^2 dx dt-\int_0^t\int_\g v^s,_3^h\cdot v^h dx_hdt \right|\le C_3\,.\label{cv17}
\end{equation}
 Using (\ref{cv16}) again and $v\in L^2(0,\infty;H^3(\of))$, we deduce from (\ref{cv17}) that
 \begin{equation}
 	\int_0^\infty 	\|v,_3^h\|^2_{L^2(\os)} dt<\infty\,.\label{cv18}
 \end{equation}
In the same way as we proved (\ref{cv1}), we infer from (\ref{cv15}) and (\ref{cv18}) that
\begin{equation*}
	\lim_{t\rightarrow\infty}	\|\eta,_3^h\|_{L^2(\os)}=0\,,
\end{equation*}
which with $\eta=0$ on $\Gamma_B$ and the fundamental theorem of calculus implies
\begin{equation}
	\lim_{t\rightarrow\infty}	\|\eta^h\|_{L^2(\os)}=0\,.\label{cv19}
\end{equation}
 Moreover, since $\|\eta\|_{L^\infty(0,\infty;H^3(\os))}<\infty$, we have by interpolation using (\ref{cv19}) that 
\begin{equation}
	\lim_{t\rightarrow\infty}	\|\eta^h\|_{H^2(\os)}=0\,.\label{cv20}
\end{equation}
From the linear wave equation in $\os$, we immediately infer from (\ref{cv20}) that
\begin{equation}
	\lim_{t\rightarrow\infty}	\|v_t^h\|_{L^2(\os)}=0\,.\label{cv21}
\end{equation}
For the velocity, (\ref{cv18})  and $v=0$ on $\Gamma_B$ imply
\begin{equation}
\int_0^\infty \|v^h\|^2_{L^2(\os)} dt<\infty\,.\label{cv22}
\end{equation}

 We prove just after that the relations (\ref{cv21}) and (\ref{cv22}) imply that
 \begin{equation}
 	\lim_{t\rightarrow\infty}	\|v^h\|_{L^2(\os)}=0\,.\label{128.11}
 \end{equation}
 To prove this to be true, let us assume by contradiction that the negation of (\ref{128.11}) is true. We then have the existence of some $\alpha >0$ such that we have a sequence $t_n\rightarrow\infty$ with $$\|v^h(t_n)\|_{L^2(\os)}>2\alpha\,.$$
  Since $$\frac{d}{dt}\|v^h\|^2_{L^2(\os)}=2\int_\os v^h\cdot v^h_t dx\,,$$
 due to (\ref{cv21}) and $N\in L^\infty(0,\infty)$, we also have the existence of $A>0$ such that for all $n$ large enough
 \begin{equation*}
 \forall t\in [t_n-A,t_n+A]\,,\ \ \|v^h(t_n)\|_{L^2(\os)}>\alpha\,,
 \end{equation*}
 which then implies in turn
 \begin{equation}
  \int_{t_n-A}^{t_n+A}\|v^h\|^2_{L^2(\os)} ds>2A\alpha^2\,.\label{128.12}
 \end{equation}
 The relations (\ref{128.12}) and (\ref{cv22}) are contradictory, which in turn implies that (\ref{128.11}) is true.
  Moreover, since $\|v\|_{L^\infty(0,\infty;H^2(\os))}<\infty$, we have by interpolation using (\ref{128.11}) that 
 \begin{equation}
 	\lim_{t\rightarrow\infty}	\|v^h\|_{H^1(\os)}=0\,.\label{cv23}
 \end{equation}
 
 \subsubsection{Convergence of one horizontal derivative of the vertical displacement and its first time derivative to zero} 
 From the continuity of stress alongside $\g$ (\ref{systemhorizontal}), we have:
 \begin{align}
 	\|\hd\eta^s,_{3}^3\|_{L^2(\g)}\lesssim& \|\nu \hd(\a_j^3\A_j^k v^f,_k^3) - \hd q \a_3^3-q \hd \a_3^3\|_{L^2(\g)}\n\\
 	\lesssim & \| v\|_{H^3(\of)} + \|\hd q\|_{L^2(\g)}+ \|q\|_{L^4(\g)} \|\hd^2\eta\|_{L^4(\g)}\n\\
 	\lesssim & \| v\|_{H^3(\of)}+ \|\hd q\|_{H^1(\of)} + \|\hd^2\tilde \eta\|_{H^{\ud}(\g)}
 	\,,\label{2310.01}
 \end{align}
 where we used (\ref{177.2bis}) to control $\|q\|_{L^4(\g)}$ independently of time.
 
 From (\ref{cv1.ter}), we know that the last term of the right hand side of (\ref{2310.01}) is square integrable over $[0,\infty)$ (for initial data satisfying Theorem \ref{theorem_main}).

Due to (\ref{0611.2}),

 
 \begin{equation}
 \int_0^\infty \|\hd q\|^2_{H^1(\of)} dt<\infty\,.\label{2310.02}
 \end{equation}
 
Using (\ref{2310.02}) in (\ref{2310.01}) then yields 
  \begin{equation}
 	\int_0^\infty 	\|\hd\eta^s,_{3}^3\|^2_{L^2(\g)} dt<\infty\,.\label{cv24}
 \end{equation}
  In a way similar as we obtained (\ref{cv9}) (due to (\ref{cv24}) playing the same role as (\ref{cv5}) did for the horizontal components), we also have the existence of $C_4>0$ finite such that
 \begin{equation}
 	\forall t\ge 0\,,\ \ \left|\int_0^t\int_\os 2|\hd\eta,_3^3|^2 dx dt-\int_0^t\int_\g \hd\eta,_3^3 \hd\eta^3 dx_hdt \right|\le C_4\,.\label{2310.03}
 \end{equation}
 From (\ref{cv1.ter}) we infer
 \begin{equation}
 \infty>\int_0^\infty \|\hd\e\|^2_{H^1(\g)} dt=\int_0^\infty \|\hd\eta\|^2_{H^1(\g)} dt\,.\label{2310.04}
 \end{equation}
 Using (\ref{cv24}) and (\ref{2310.04}) in (\ref{2310.03}), we obtain
 \begin{equation}
 	\int_0^\infty 	\|\hd\eta,_{3}^3\|^2_{L^2(\os)} dt<\infty\,.\label{cv25}
 \end{equation} 
Since $\hd\eta^3=0$ on $\Gamma_B$, we infer from (\ref{cv25}) and the fundamental theorem of calculus that
\begin{equation}
	\int_0^\infty 	\|\hd\eta^3\|^2_{L^2(\os)} dt<\infty\,.\label{cv26}
\end{equation}
Given that
\begin{equation}
\left\|\ddt	\|\hd\eta^3\|^2_{L^2(\os)}\right\|_{L^\infty(0,\infty)} =2 \left\|\int_\os \hd\eta^3 \hd v^3 dx\right\|_{L^\infty(0,\infty)}\lesssim N(t)\le C<\infty\,,\label{cv27}
\end{equation}
we infer from (\ref{cv26}), (\ref{cv27}) (in the same way as \ref{cv23}) that
\begin{equation}
	\lim_{t\rightarrow\infty}	\|\hd\eta^3\|^2_{L^2(\os)} =0\,.\label{cv28}
\end{equation}
 Moreover, since $\|\hd\eta\|_{L^\infty(0,\infty;H^2(\os))}<\infty$, we have by interpolation using (\ref{cv28}) that 
\begin{equation}
	\lim_{t\rightarrow\infty}	\|\hd\eta^3\|_{H^1(\os)}=0\,.\label{cv29}
\end{equation}


We now prove that
\begin{equation}
	\lim_{t\rightarrow\infty}	\|\hd v^3\|_{L^2(\os)}=0\,.\label{cv31}
\end{equation}
In order to do so, we proceed by contradiction, and assume that there exists $M>0$ finite such that
\begin{equation}
	\limsup_{t\rightarrow\infty}	\|\hd v^3\|^2_{L^2(\os)}=4M>0\,.\label{cv32}
\end{equation}
Therefore, we have a sequence $t_n\rightarrow\infty$ such that 
\begin{equation}
\forall n\in \mathbb{N}\,,\ \ \|\hd v^3\|^2_{L^2(\os)}(t_n)\ge 2 M\,.\label{cv33}
\end{equation}
Since there exists $C_5>0$ finite such that
\begin{equation}
\forall t\ge 0\,,\ 	\left|\ddt\|\hd v^3\|^2_{L^2(\os)}\right|=2\left|\int_\os \hd v^3 \hd v_t^3 dx\right| \le C_5\,,\label{cv34}
\end{equation}
we infer from (\ref{cv33}), (\ref{cv34}) and the fundamental theorem of calculus that
\begin{equation}
	\forall t\in [t_n,t_n+\frac{M}{C_5}]\,,\ \ \|\hd v^3\|^2_{L^2(\os)}(t)\ge  M\,.\label{cv35}
\end{equation}

We now take one horizontal derivative of the vertical component of the wave equation (which makes the gravity term disappear), and multiply by $\hd \eta^3$:
\begin{align}
	0=&\int_{t_n}^{t_n+\frac{M}{C_5}}\int_\os \hd v_t^3 \hd\eta^3 dx dt+\lambda\int_{t_n}^{t_n+\frac{M}{C_5}}\int_\os |\hd\nabla\eta^3|^2 dxdt -\lambda\int_{t_n}^{t_n+\frac{M}{C_5}}\int_\g \hd\eta,_{3}^3 \hd\eta^3 dx_h dt\n\\
	=&-\int_{t_n}^{t_n+\frac{M}{C_5}}\int_\os |\hd v^3|^2dx dt+\left[\int_\os \hd v^3\hd\eta^3 dx\right]_{t_n}^{t_n+\frac{M}{C_5}}+\lambda\int_{t_n}^{t_n+\frac{M}{C_5}}\int_\os |\nabla\hd\eta^3|^2 dxdt\n\\ &-\lambda\int_{t_n}^{t_n+\frac{M}{C_5}}\int_\g \hd\eta,_{3}^3 \hd\eta^3 dx_hdt\,.\label{cv36}
\end{align}

Due to (\ref{cv29}) and (\ref{cv32}), we have that the second term on the right hand side of (\ref{cv36}) converges to zero as $n\rightarrow\infty$. The same holds true for the third term as well due to (\ref{cv29}). Due to (\ref{2310.04}) and (\ref{cv24}), we have the same convergence for the last term on the right hand side of (\ref{cv36}) as well. On the other hand, for the first term on the right hand side of (\ref{cv36}), we have thanks to (\ref{cv35}):
$$0<M \frac{M}{C_5 }\le \int_{t_n}^{t_n+\frac{M}{C_5}}\int_\os |\hd v^3|^2dx dt \rightarrow 0\ \ \text{as}\ n\rightarrow\infty\,.$$
This is clearly not possible, which leads us to reject our assumption (\ref{cv32}) that led to this contradiction. This in turn establishes that (\ref{cv31}) is true.
 
\subsubsection{Convergence of the solution towards a flat interface solution as $t\rightarrow\infty$}
 
 We now prove that there exists $(\alpha_0,\alpha_1)\in H_0^1(0,h_s)\times H_0^1(0,h_s)$ such that if we denote by $\alpha$ the
 solution of the $1-$d wave equation (\ref{wave1d})
 then,
 \begin{equation}
 \lim_{t\rightarrow\infty} (\|\alpha(\cdot,t)-\eta^3(\cdot,\cdot,t)\|_{H^1(\os)} + \|\alpha_t(\cdot,t)-v^3(\cdot,\cdot,t)\|_{L^2(\os)})=0\,.\label{39.1}
 \end{equation}
 
 The proof will be split into five steps.
 
  In Step 1, we define elliptic operators $\Lambda_i$ ($i=0,1$) in $\os$ associated with a Dirichlet boundary condition set to zero on the reference interface $\g$. We also define a sequence $\alpha^n$ of one dimensional wave problems (with Dirichlet boundary condition set to zero), with data specified at time $n$ as being the horizontal average of the solution of (\ref{ale}) at time $n$. 
  
  In Step 2, we show the existence of a weakly convergent subsequence of the solutions of the one dimensional problems defined in Step 1. This provides us with $\alpha_0$ and $\alpha_1$.
  
  In Step 3, we show that $\Lambda_i$ is converging toward the identity map in large time when applied to the solution of the parabolic-hyperbolic interaction problem (\ref{ale}).
  
  In Step 4, we establish that in large time, the difference between the solution of (\ref{ale}) and its horizontal average converges towards $0$.
  
  In Step 5, we conclude by estimates on the wave equation satisfied by the difference between $\alpha^n$ and the horizontal average of the solution of (\ref{ale}).
 
 \begin{remark}
 Note that in (\ref{39.1}), we view $\alpha(x_3,t)$ as being a constant function in the horizontal variable $x_h=(x_1,x_2)$. This convention will be used later on as well.
 \end{remark}
 
 \noindent{\bf Step 1. Definition of the horizontal average, of an elliptic problem, and of a wave system with data prescribed at time $n$.}
 
 \begin{definition}
 It will be convenient to introduce in this section a notation for the horizontal average of a function $F$, that we denote as
 \begin{equation*}
 m(F)(x_3,t)=\frac{1}{L^2} \int_{[0,L]^2} F(x_h,x_3,t)\ dx_h\,.
 \end{equation*}
\end{definition}

\begin{definition}
 We next introduce for any smooth enough $f$, and for each $i=0,1$ the associated $\Lambda_i f$ defined as the periodic solution (in the canonical horizontal directions) of the elliptic system:
 \begin{subequations}
 	\label{h10}
 	\begin{alignat}{2}
 		\Delta \Lambda_i f &=\Delta f \ \  \ &&\text{in} \ \ \os \,,\\
 		\Lambda_i f &= (h_e-h_s)\delta_i^0 \ \  &&\text{on} \ \ \g \,,\label{h10.b}\\
 			\Lambda_i f &= 0 \ \  &&\text{on} \ \ \g\cup\Gamma_B \,,\label{h10.b1}
 	\end{alignat}
 \end{subequations}
\end{definition} 
 
  By standard variational regularity, we have that 
 \begin{equation}
 	\label{119.1}
 	\|\Lambda_i f\|_{H^1(\os)}\lesssim \|f\|_{H^1(\os)}+|h_e-h_s| \delta_i^o\,,
 \end{equation}
 as well as
 \begin{equation}
 	\label{100825.1}
 	\|\Lambda_i f\|_{H^2(\os)}\lesssim \|\Delta f\|_{L^2(\os)}+|h_e-h_s|\delta_i^0\,.
 \end{equation}

\begin{definition}
 Next, we define for each $n\in\mathbb{N}$, $\alpha^n$ as the solution of the one dimensional wave equation with data prescribed at time $n$:
\begin{subequations}
 	\label{waven1d}
 	\begin{alignat}{2}
 		\alpha_{tt}^n-\lambda\alpha,_{33}^n &=-g \ \  \ &&\text{in} \ \ (0,h_s)\times [0,\infty) \,,\\
 		\alpha^n &= h_e-h_s \ \  &&\text{on} \ \ \{h_s\} \times [0,\infty)\,,\\
 		\alpha^n &= 0 \ \  &&\text{on} \ \ \{0\} \times [0,\infty)\,,\\
 		\alpha^n(x_3,n) &=m(\Lambda_0\eta^3)(x_3,n))=\frac{1}{L^2}\int_{[0,L]^2} \Lambda_0 \eta^3 (x_h,x_3,n) dx_h \ \   &&\text{on} \ \ [0,h_s]\,,\\
 		\alpha^n_t(x_3,n) &=m(\Lambda_1 v^3)(x_3,n)=\frac{1}{L^2}\int_{[0,L]^2} \Lambda_1 v^3 (x_h,x_3,n) dx_h \ \   &&\text{on} \ \ [0,h_s]\,.
 	\end{alignat}
 \end{subequations}
\end{definition}


Since the wave equation is time revertible, the data being prescribed at time $n$ still allows for the solution being defined for all time.

Due to (\ref{h10.b}), (\ref{h10.b1}) we have that $\alpha^n(0,n)=0$ and $\alpha_n(h_s,n)=h_e-h_s$, which ensures that the solution of (\ref{waven1d}) has the regularity
$(\alpha^n,\alpha^n_t)\in L^\infty (0,\infty; H^1_0(0,h_s))\times L^\infty (0,\infty; L^2(0,h_s))$.

Similarly, due to (\ref{h10.b}), (\ref{h10.b1}) we have that $\alpha_t^n(\cdot,n)=0$ on $\{0\}\cup\{h_s\}$, which ensures that the solution of the one time differentiated version of (\ref{waven1d}) has the regularity
$(\alpha^n_t,\alpha^n_{tt})\in L^\infty (0,\infty; H^1_0(0,h_s))\times L^\infty (0,\infty; L^2(0,h_s))$, with the estimate a.e on $[0,\infty)$:
\begin{align}
\|\alpha^n_{tt}\|_{L^2(0,h_s)}^2+ \lambda\|\alpha^n_{t},_3\|_{L^2(0,h_s)}^2= & \|\alpha^n_{tt}(\cdot,n)\|_{L^2(0,h_s)}^2+ \lambda\|\alpha^n_{t},_3(\cdot,n)\|_{L^2(0,h_s)}^2\n\\
= & \|\lambda m((\Lambda_0 \eta^3),_{33})(\cdot,n)-g\|_{L^2(0,h_s)}^2+\lambda\|m((\Lambda_1 v^3),_3)(\cdot,n)\|_{L^2(0,h_s)}^2\n\\
\lesssim & \|(\Lambda_0 \eta^3),_{33}(\cdot,n)\|_{L^2(\os)}^2+g^2+\|(\Lambda_1 v^3),_3(\cdot,n)\|_{L^2(\os)}^2\n\\
\lesssim & \|\Lambda_0 \eta^3(\cdot,n)\|_{H^2(\os)}^2+g^2+\|\Lambda_1 v^3(\cdot,n)\|_{H^1(\os)}^2\n\\
\lesssim & \|\Delta\eta^3(\cdot,n)\|_{L^2(\os)}^2+|h_e-h_s|^2+g^2+\|v^3(\cdot,n)\|_{H^1(\os)}^2\n\\
\lesssim & \| v_t(\cdot,n)\|_{L^2(\os)}^2+|h_e-h_s|^2+g^2+\|v^3(\cdot,n)\|_{H^1(\os)}^2
\,,\label{119.2}
\end{align}
where we used the elliptic regularity (\ref{119.1}), (\ref{100825.1}) to obtain the inequlaity before (\ref{119.2}) from the previous inequality. Note also that in the second equality above we just used the fact $m$ commutes with the vertical derivative, and in the inequality below the fact that the $L^2$ norm of the average of a function is controlled by the $L^2$ norm of the function by Cauchy-Schwarz.
Remember also that gravity disappears in the first time-differentiated wave equation. Using our global in time estimate, we then infer from (\ref{119.2}) that
\begin{equation}
\|\alpha^n_t\|_{L^\infty(0,\infty;H^1_0(0,h_s))}^2+ \|\alpha^n_{tt}\|_{L^\infty(0,\infty;L^2(0,h_s))}^2\lesssim \epsilon_0+g^2+\underbrace{|h_e-h_s|^2}_{\lesssim \epsilon_0}\,.\label{119.3}
\end{equation}
Similarly,
 \begin{equation}
 	\|\alpha^n\|_{L^\infty(0,\infty;H^1(0,h_s))}^2+ \|\alpha^n_{t}\|_{L^\infty(0,\infty;L^2(0,h_s))}^2\lesssim \epsilon_0+g^2+\underbrace{|h_e-h_s|^2}_{\lesssim \epsilon_0}\,,\label{119.4}
 \end{equation}
 due to the presence of gravity (either zero or non zero) in (\ref{waven1d}). Due to $\alpha$ being solution of the wave equation, we infer from (\ref{119.3}) that
 \begin{equation}
 	 \|\alpha^n,_{33}\|_{L^\infty(0,\infty;L^2(0,h_s))}^2\lesssim \epsilon_0+g^2\,.\label{1611.4}
 \end{equation}
 
 \begin{remark}
 Note that the compatibility conditions on $\alpha(\cdot,n)$ and $\alpha_t(\cdot,n)$ on $\{0\}\cup\{h_s\}$ being satisfied at time $n$ ensure that they are satisfied for all time, and in particular at time $0$.
  \end{remark}
 
\noindent{\bf Step 2. Convergence of a subsequence solution of (\ref{waven1d}) towards a solution of (\ref{wave1d}).}

By (\ref{119.3}), (\ref{119.4}) and (\ref{1611.4}), we have that $e^{-t} \alpha^n$ is bounded in $H^2((0,h_s)\times(0,\infty))$. By weak convergence in $H^2((0,h_s)\times(0,\infty))$, we have the existence of a strictly increasing mapping $\sigma$ from $\mathbb{N}$
into itself such that $e^{-t}\alpha^{\sigma(n)}$ converges weakly in $H^2((0,h_s)\times(0,\infty))$ towards some limit $e^{-t}\alpha\in H^2((0,h_s)\times(0,\infty))$. By compactness of the trace theorem from $H^1((0,h_s)\times (0,\infty))$ into $L^2((0,h_s)\times\{0\})$, we then infer that
$$\|\alpha^{\sigma(n)}-\alpha\|^2_{L^2((0,h_s)\times\{0\})}+ \|\alpha_t^{\sigma(n)}-\alpha_t\|^2_{L^2((0,h_s)\times\{0\})}
+ \|\alpha,_3^{\sigma(n)}-\alpha,_3\|^2_{L^2((0,h_s)\times\{0\})}\rightarrow 0\ \text{as}\ n\rightarrow\infty\,.$$
 
 From this, we then have the existence of $(\alpha_0,\alpha_1)\in H^1(0,h_s)\times L^2(0,h_s)$ such that
 \begin{subequations}
 	\label{79.4}
 	\begin{align}
 		\alpha_t^{\sigma(n)}(0)&\rightarrow \alpha_1=\alpha_t(0)\ \text{in}\ L^2(0,h_s)\,,\\
 		\alpha^{\sigma(n)}(0)&\rightarrow \alpha_0=\alpha(0)\ \text{in}\ H^1(0, h_s)\,.
 	\end{align}
 \end{subequations}
 
 Due to the boundary conditions on (\ref{waven1d}), we also have $\alpha_0(0)=0$ and $\alpha(h_s)=h_e-h_s$.
 
Moreover, due to  (\ref{119.3}), we also have $\alpha_1\ \text{in}\ H^1_0(0, h_s)$ with the weak convergence
\begin{equation*}
	\alpha_t^{\sigma(n)}(0)\rightharpoonup \alpha_1\ \text{in}\ H^1_0(0,h_s)\,.
\end{equation*}
 
 We also have similar convergence properties at any time $t>0$ (by replacing $0$ by $t$ above):
 	\begin{align*}
 		e^{-t}\alpha_t^{\sigma(n)}(t)&\rightarrow e^{-t}\alpha_t(t)\ \text{in}\ L^2(0,h_s)\,,\\
 		e^{-t}\alpha^{\sigma(n)}(t)&\rightarrow e^{-t}\alpha(t)\ \text{in}\ H^1(0, h_s)\,.
 	\end{align*} 
 which implies
  \begin{subequations}
 	\label{109.2}
 	\begin{align}
 		\alpha_t^{\sigma(n)}(t)&\rightarrow \alpha_t(t)\ \text{in}\ L^2(0,h_s)\,,\\
 		\alpha^{\sigma(n)}(t)&\rightarrow \alpha(t)\ \text{in}\ H^1(0, h_s)\,.
 	\end{align}
 \end{subequations}
 
 Due to the bounds (\ref{119.3}) and (\ref{1611.4}), we also have that $\alpha_{tt}$ and $\alpha,_{33}$ are in $L^\infty(0,\infty;L^2(0,h_s))$.

Due to each $\alpha^n$ satisfying the system (\ref{waven1d}), we have as a result of these convergences that $\alpha$ is solution of the wave system (\ref{wave1d}) with initial data $(\alpha_0,\alpha_1)$.

\noindent{\bf Step 3. Convergence of $(\eta_3-\Lambda_0\eta_3,v_3-\Lambda_0 v_3)$ towards $0$.} 
 
By definition of our Dirichlet problem (\ref{h10}), we have
\begin{align}
\|\nabla (\Lambda_1 v_3-v_3)\|^2_{L^2(\os)}=&\int_{\g\cup\Gamma_B} (\Lambda_1 v_3-v_3),_3 (\Lambda_1 v_3-v_3) N_3 dx_h\n\\
 =&\int_{\g} (\Lambda_1 v_3-v_3),_3 (\Lambda_1 v_3-v_3) dx_h\n\\
 =&-\int_{\g} (\Lambda_1 v_3-v_3),_3 v_3 dx_h\n\\
 \le\ & \|(\Lambda_1 v_3-v_3),_3\|_{L^2(\g)} \|v_3\|_{L^2(\g)}\,.\label{59.1}
\end{align}
 By $H^{\frac 3 2}(\os)$ regularity on the elliptic problem (\ref{h10}), we infer successively from (\ref{59.1}) and Theorem \ref{theorem_main} that
 \begin{equation}
 	\|\nabla (\Lambda\-1 v_3-v_3)\|^2_{L^2(\os)}
 	\lesssim  \|v_3\|_{H^2(\os)} \|v_3\|_{L^2(\g)}
 	\lesssim  \|v_3\|_{L^2(\g)}
 	\,.\label{59.2}
 \end{equation}
 Due to the velocity in the fluid phase converging to $0$ by (\ref{cv1}), we infer from (\ref{59.2})
 \begin{equation}
 \lim_{t\rightarrow\infty} \|\Lambda_1 v_3-v_3\|_{H^1(\os)}=0\,.\label{59.3}
 \end{equation}
We obtain similarly as (\ref{59.1}) that
 $$\|\nabla (\Lambda_0 \eta_3-\eta_3)\|^2_{L^2(\os)}\lesssim \|(\Lambda_0 \eta_3-\eta_3),_3\|_{L^2(\g)} \|\Lambda_0\eta^3-\eta_3\|_{L^2(\g)}= \|(\Lambda_0 \eta_3-\eta_3),_3\|_{L^2(\g)} \|h_e-h_s-\eta_3\|_{L^2(\g)}\,,$$
 which allows us to infer from (\ref{cv4}) that
  \begin{equation}
 	\lim_{t\rightarrow\infty} \|\Lambda_0 \eta_3-\eta_3\|_{H^1(\os)}=0\,.\label{59.4}
 \end{equation}
 
 \noindent{\bf Step 4. Convergence of $(\eta_3-m(\eta_3),v_3- m(v_3))$ towards $0$.}

 Due to Poincar\'e-Wirtinger in $[0,L]^2$, we have for each $x_3\in (0,h_s)$
 \begin{align*}
 \|(m(\eta,_3^3)-\eta,_3^3)(\cdot,x_3,t)\|^2_{L^2((0,L)^2)}\lesssim \|\hd\eta,_{3}^3(\cdot,x_3,t)\|^2_{L^2((0,L)^2)}\,,
 \end{align*}
 which by vertical integration yields
 \begin{align}
 	\|(m(\eta^3)-\eta^3),_3(\cdot,t)\|^2_{L^2(\os)}\lesssim \|\hd\eta,_{3}^3(\cdot,t)\|^2_{L^2(\os)}\,.\label{79.5}
 \end{align}
 On the other hand, since $\hd m(\eta^3)=0$,
 \begin{align}
 	\|\hd(m(\eta^3)-\eta^3)(\cdot,t)\|^2_{L^2(\os)}= \|\hd\eta^3(\cdot,t)\|^2_{L^2(\os)}\,.\label{79.6}
 \end{align}
 Therefore, (\ref{79.5}), (\ref{79.6}) and the convergence to $0$ of $\hd\eta^3$ in $H^1(\os )$ by (\ref{cv29}) yield
 \begin{equation*}
 \lim_{t\rightarrow\infty} \|\nabla(m(\eta^3)-\eta^3)\|_{L^2(\os)}=0\,,
 \end{equation*}
 which by Poincar\'e's inequality (due to $\eta_3=0=m(\eta^3)$ on $\Gamma_B$) implies
 \begin{equation}
 	\lim_{t\rightarrow\infty} \|m(\eta^3)-\eta^3\|_{H^1(\os)}=0\,,\label{79.7}
 \end{equation}
 Similarly, due to (\ref{cv31}) we also have in $L^2(\os)$ (without taking a vertical derivative)
 \begin{equation}
 	\lim_{t\rightarrow\infty} \|m(v^3)-v^3\|_{L^2(\os)}\lesssim 	\lim_{t\rightarrow\infty} \|\hd v^3\|_{L^2(\os)}= 0\,.\label{79.8}
 \end{equation}
 
 \noindent{\bf Step 5. Study of the difference between $\alpha^n$ and $m(\eta^3)$, and between  $\alpha^n_t$ and $m(v^3)$.} 
  To shorten notations, we denote
 \begin{equation*}
 \delta(t)=\|\alpha^n_t-m(v^3)\|^2_{L^2(0,h_s)}(t)+ \lambda\|(\alpha^n-m(\eta^3)),_3\|^2_{L^2(0,h_s)}(t)\,.
 \end{equation*}
 We observe that (using the value at time $n$ in (\ref{waven1d}), and the fact the $L^2$ norm of the horizontal average of a function is controlled by the $L^2$ norm of the function):
 \begin{align}
 	\delta(n)=& \|\alpha^n_t-m(v^3)\|^2_{L^2(0,h_s)}(n)+ \lambda\|(\alpha^n-m(\eta^3)),_3\|^2_{L^2(0,h_s)}(n)\n\\
 	=& \|m(\Lambda_1 v^3-v_3)\|^2_{L^2(0,h_s)}(n)+ \lambda\|m( (\Lambda_0 \eta^3-\eta_3),_3)\|^2_{L^2(0,h_s)}(n)\n\\
 	\lesssim & \|\Lambda_1 v^3-v_3\|^2_{L^2(\os)}(n)+ \lambda\| (\Lambda_0 \eta^3-\eta_3),_3\|^2_{L^2(\os)}(n)\n\\
 	\rightarrow & 0\ \text{as}\ n\rightarrow\infty
 	\,,\label{2410.02}
 \end{align}
  due to (\ref{59.3}) and (\ref{59.4}).
  
 The basic energy estimate on the system (note the gravity terms cancel each other when forming the difference):
 \begin{subequations} 	\label{wave1d.bis}
 	\begin{alignat}{2}
 		(\alpha^n-m(\eta^3))_{tt}-\lambda(\alpha^n-m(\eta^3)),_{33} &=0 \ \  \ &&\text{in} \ \ (0,h_s)\times [0,\infty) \,,\\
 		\alpha^n-m(\eta^3) &= -m(\eta^3) \ \  &&\text{on} \ \ \{0,h_s\}\times [0,\infty) \,,
 	\end{alignat}
 \end{subequations}
yields 
 \begin{align}
 \delta(t)= & \delta(n)+ 2\lambda \int_n^t  (\alpha^n-m(\eta^3)),_{3}(h_s,s) (\alpha^n_t-m (v^3))(h_s,s) ds\n\\
 = & \delta(n)-2\lambda \int_n^t  (\alpha^n-m(\eta^3)),_{3}(h_s,s) m(v^3)(h_s,s)  ds\n\\
  = & \delta(n) -2\lambda \int_n^t  (\alpha^n-m(\eta^3)),_{3}(h_s,s) m(v^3)(h_s,s) \frac{x_3}{h_s}(h_s,s)  ds
 \,,\label{89.1}
 \end{align}
 where we simply used $\frac{x_3}{h_s}=1$ on $x_3=h_s$ in (\ref{89.1}).
 
 By the fundamental theorem of calculus for the function $x_3\rightarrow (\alpha^n-m(\eta^3)),_{3}(x_3,t) m(v^3)(h_s,t) \frac{x_3}{h_s}$ (which is zero for $x_3=0$), we infer from (\ref{89.1}) that
 \begin{align}
 	\delta(t)
 	&=  \delta(n) -2\lambda \int_n^t \int_0^{h_s} (\alpha^n-m(\eta^3)),_{33}  m(v^3)(h_s,s) \frac{x_3}{h_s} + (\alpha^n-m(\eta^3)),_{3} m(v^3)(h_s,s) \frac{1}{h_s} dx_3  ds\n\\
 	&= \delta(n)-2 \int_n^t \int_0^{h_s} (\alpha^n-m(\eta^3)),_{tt} m(v^3)(h_s,s) \frac{x_3}{h_s} + \lambda (\alpha^n-m(\eta^3)),_{3} m( v^3)(h_s,s) \frac{1}{h_s} dx_3  ds\n\\
 	&\lesssim \delta(n)+\left|\int_n^t |m(v^3)(h_s,s)| ds\right|\,,
 	\label{89.2}
 \end{align}
 where we used (\ref{119.3}) and (\ref{119.4}), as well as the previously established control for all time of $v^3_t$ and $\eta,_3^3$ in $L^2(\os)$, in order to obtain (\ref{89.2}) from the previous line.

 We now explain why the integral on the right hand side of (\ref{89.2}) converges as $t\rightarrow\infty$, despite not being the integral of a square. Using the divergence free condition, and $v_3=0$ on the fixed $\Gamma_{top}$ we have by integration by parts in the fluid reference domain:
 \begin{align*}
 -\int_\g v^3 dx_h=& \int_\g v^3 N^3 dx_h\n\\
 =& \int_\g v^i N^i dx_h\n\\
  =& -\int_{\Gamma_{top}} \u{v^i}_{=0\ \text{on}\ \Gamma_{top}} N^i dx_h + \int_\of \div v\   dx\n\\
  =& \int_\of (\delta_i^j-\u{\a_i^j)v,_j^i}_{=0}   dx\n\\
  =& -\int_\of (\delta_i^j-\a_i^j),_j v^i   dx+ \int_{\p\of} (\delta_i^j-\a_i^j) \u{v^i}_{=0\ \text{on}\ \Gamma_{top}} N_j dx_h\n\\
  =& \int_\of \u{\a_i^j,_j}_{=0\ \text{by (\ref{piola})}} v^i   dx- \int_{\g} (\delta_i^3-\a_i^3) v^i  dx_h\,.
 \end{align*}
Since $(\a_i^3)_{i=1}^3=(\text{Id}+\e),_1\times(\text{Id}+\e),_2$, we infer from this equality and Cauchy-Schwarz that
 \begin{equation*}
 |m(v^3)(h_s,t)|\lesssim \|\hd\tilde \eta\|_{L^2(\g)}\|v\|_{H^1(\of)}\,,
 \end{equation*}
 which then implies with (\ref{cv0}) and (\ref{cv1.ter}) that
 \begin{equation}
 \int_0^\infty	|m(v^3)(h_s,s)| ds<\infty\,.\label{109.1}
 \end{equation}
 Now, let $\varepsilon>0$ be fixed. Using (\ref{109.1}), we have the existence of $T_\varepsilon>0$ such that
 \begin{equation}
  \int_{T_\epsilon}^\infty	|m(v^3)(h_s,s)| ds<\varepsilon\,.\label{1611.1}
 \end{equation}
 From (\ref{1611.1}), we infer that
 \begin{equation}
 	\forall t\ge T_\varepsilon\,,	\forall n\ge T_\varepsilon\,,\ \left|\int_n^t	|m(v^3)(h_s,s)\right| ds |<\varepsilon\,.\label{1611.2}
 \end{equation}
 Due to the convergence (\ref{2410.02}), we can also have $N_{\varepsilon}$ large enough so that
 \begin{equation}
 	\forall n\ge N_{\varepsilon}\,, \ \delta(n)<\varepsilon\,.\label{1611.3}
 \end{equation}

 Using (\ref{1611.2}) and (\ref{1611.3}) in (\ref{89.2}) yields
 \begin{equation*}
 \forall t\ge T_\varepsilon\,,\  \forall n\ge \max(T_\varepsilon,N_\varepsilon)\,,	\|(\alpha^n_t-m( v^3))(\cdot,t)\|^2_{L^2(0,h_s)}+ \lambda\|(\alpha^n-m( \eta^3)),_3(\cdot,t))\|^2_{L^2(0,h_s)}\lesssim \varepsilon\,.
 \end{equation*}
 As $\sigma(n)\ge n$ for an increasing injective integer map, this estimate is again true for $\sigma(n)$:
 \begin{equation}
 	\forall t\ge T_\varepsilon\,,\  \forall n\ge \max(T_\varepsilon,N_\varepsilon)\,,	\|(\alpha^{\sigma(n)}_t-m( v^3))(\cdot,t)\|^2_{L^2(0,h_s)}+ \lambda\|(\alpha^{\sigma(n)}-m( \eta^3)),_3(\cdot,t))\|^2_{L^2(0,h_s)}\lesssim \varepsilon\,.\label{109.3}
 \end{equation}

 Using the convergences (\ref{109.2}), and the estimate (\ref{109.3}), we infer:
 \begin{equation}
 	\forall t\ge T_\varepsilon\,,\ \ 	\|(\alpha_t- m(v^3))(\cdot,t)\|^2_{L^2(0,h_s)}+ \lambda\|(\alpha-m(\eta^3)),_3(\cdot,t)\|^2_{L^2(0,h_s)}\lesssim \varepsilon\,,
 \end{equation}
 which, with the convergences (\ref{79.7}) and (\ref{79.8}) shows that
 \begin{equation}
 \lim_{t\rightarrow\infty} \|(\alpha_t- v^3)(\cdot,t)\|^2_{L^2(\os)}+ \lambda\|(\alpha-\eta^3),_3(\cdot,t)\|^2_{L^2(\os)}=0\,.\label{1611.6}
 \end{equation} 
 Given the convergence of $\hd\eta^3$ established in (\ref{cv29}), the fact that $\hd \alpha=0$, and $\alpha-\eta^3=0$ on $\Gamma_B$, (\ref{1611.6}) allows us to conclude that the convergence (\ref{39.1}) holds.
 
 This finishes the proof of Theorem 2.\qed
 
 \begin{remark}
 The fact the convergence (\ref{39.1}) holds implies the couple $(\alpha_0,\alpha_1)$ is unique and so the convergences (\ref{79.4}) and (\ref{109.2}) hold for the whole sequence $\alpha^n$.
 \end{remark}

  \end{document}